\documentclass[11pt]{article}

\usepackage{epsfig}
\usepackage{graphicx}
\usepackage{amsbsy}
\usepackage{amsmath}
\usepackage{amsfonts}
\usepackage{amssymb}
\usepackage{textcomp}
\usepackage{hyperref}
\usepackage{aliascnt}

\newcommand{\mcm}[3]{\newcommand{#1}[#2]{{\ensuremath{#3}}}} 

\mcm{\tuple}{1}{\langle #1 \rangle}
\mcm{\name}{1}{\ulcorner #1 \urcorner}
\mcm{\Nbb}{0}{\mathbb{N}}
\mcm{\Zbb}{0}{\mathbb{Z}}
\mcm{\Rbb}{0}{\mathbb{R}}
\mcm{\Cbb}{0}{\mathbb{C}}
\mcm{\Qbb}{0}{\mathbb{Q}}
\mcm{\Acal}{0}{\cal A}
\mcm{\Bcal}{0}{\cal B}
\mcm{\Ccal}{0}{\cal C}
\mcm{\Dcal}{0}{\cal D}
\mcm{\Ecal}{0}{\cal E}
\mcm{\Fcal}{0}{\cal F}
\mcm{\Gcal}{0}{\cal G}
\mcm{\Hcal}{0}{\cal H}
\mcm{\Ical}{0}{\cal I}
\mcm{\Jcal}{0}{\cal J}
\mcm{\Kcal}{0}{\cal K}
\mcm{\Lcal}{0}{\cal L}
\mcm{\Mcal}{0}{\cal M}
\mcm{\Ncal}{0}{\cal N}
\mcm{\Ocal}{0}{{\cal O}}
\mcm{\Pcal}{0}{{\cal P}}
\mcm{\Qcal}{0}{{\cal Q}}
\mcm{\Rcal}{0}{{\cal R}}
\mcm{\Scal}{0}{{\cal S}}
\mcm{\Tcal}{0}{{\cal T}}
\mcm{\Ucal}{0}{{\cal U}}
\mcm{\Vcal}{0}{{\cal V}}
\mcm{\Wcal}{0}{{\cal W}}
\mcm{\Xcal}{0}{{\cal X}}
\mcm{\Ycal}{0}{{\cal Y}}
\mcm{\Mfrak}{0}{\mathfrak M}

\mcm{\restric}{0}{\upharpoonright}
\mcm{\upset}{0}{\uparrow}
\mcm{\onto}{0}{\twoheadrightarrow}
\mcm{\smallNbb}{0}{{\small \mathbb{N}}}
\DeclareMathOperator{\preop}{op}
\mcm{\op}{0}{^{\preop}}

\newcommand{\se}{\subseteq}
{\begin{array}{c}
\setlength{\unitlength}{1em}}%
{\end{array}}

\usepackage{amsthm}

\newcommand{\theoremize}[2]{\newaliascnt{#1}{thm} \newtheorem{#1}[#1]{#2} \aliascntresetthe{#1}}

\theoremstyle{plain}
\newtheorem{thm}{Theorem}[section]
\theoremize{lem}{Lemma}
\theoremize{skolem}{Skolem}
\theoremize{fact}{Fact}
\theoremize{sublem}{Sublemma}
\theoremize{claim}{Claim}
\theoremize{obs}{Observation}
\theoremize{prop}{Proposition}
\theoremize{cor}{Corollary}
\theoremize{que}{Question}
\theoremize{oque}{Open Question}
\theoremize{con}{Conjecture}

\theoremstyle{definition}
\theoremize{dfn}{Definition}
\theoremize{rem}{Remark}
\theoremize{eg}{Example}
\theoremize{exercise}{Exercise}
\theoremstyle{plain}

\usepackage{verbatim}
\usepackage{enumerate}
\usepackage[all]{xy}

\usepackage{subfig}

\title{Embedding simply connected \\ 2-complexes in 3-space\\ \Large IV. Dual matroids}

\author{Johannes Carmesin
\medskip 
\\
  {University of Cambridge}
}

\DeclareMathOperator{\Sbb}{\mathbb{S}}

\newcommand{\sm}{\setminus}

\newcommand{\Sthree}{$\Sbb^3$}

\mcm{\Fbb}{0}{\mathbb{F}}

\usepackage{xr}
\begin{document}

\maketitle
\begin{abstract}
We introduce dual matroids of 2-dimensional simplicial complexes. 
Under certain necessary conditions, duals matroids are used to characterise embeddability in 
3-space in a way analogous to Whitney's planarity criterion. 

We further use dual matroids to extend a 3-dimensional analogue of Kuratowski's theorem to the 
class of 2-dimensional simplicial complexes obtained from simply connected ones by identifying 
vertices or 
edges. 
\end{abstract}
\section{Introduction}
A well-known characterisation of planarity of graphs is Whitney's theorem from 1932. 
It states that a graph can be embedded in the plane if and only if its dual 
matroid is \emph{graphic} (that is, it is the cycle matroid of a graph) \cite{Whitney32}. 

In this paper we define dual matroids of (2-dimensional) simplicial complexes. We prove 
under certain necessary assumptions 
an analogue of Whitney's characterisation for embedding simplicial complexes in 3-space. More 
precisely, 
under these assumptions a simplicial complex can be embedded in 3-space if and only if its dual 
matroid 
is graphic.

\vspace{.3 cm}

Our definition of dual matroid is inspired by the following fact. 
\begin{thm}\label{3dual_matroid}
Let $C$ be a directed 2-dimensional simplicial complex embedded into \Sthree. Then the edge/face 
incidence matrix of $C$ 
represents over the integers\footnote{See \autoref{sec:dual_matroids} for a definition. } a 
matroid $M$ which is equal to the cycle matroid of the dual graph of the embedding.  
\end{thm}
Indeed, we define\footnote{The choice of $\Fbb_3$ is a bit 
arbitrary. Indeed any other field $\Fbb_p$ with $p$ a prime different from $2$ works. }. the 
\emph{dual 
matroid} of a simplicial complex $C$ to be the matroid represented by the edge/face incidence matrix 
of $C$ 
over 
the finite field $\Fbb_3$. 

\vspace{.3 cm}

Although the cone over $K_5$ does not embed in 3-space\footnote{See for example \cite{3space1}. }, 
its dual matroid just consists of a bunch 
of loops, and thus is graphic. In 
order to exclude examples like the cone over $K_5$ we restrict our attention to simplicial complexes 
$C$ 
whose dual matroid captures the local structure at all vertices of $C$. We call such dual 
matroids \emph{local}, see \autoref{secW} for a precise definition. 
Examples of simplicial complex whose 
dual matroid is local are those where every edge is incident with precisely three faces and the 
dual matroid has no loops. Another example is the 
3-dimensional grid whose faces are the 4-cycles. 

Furthermore matroids (of graphs and also of 
simplicial complexes) do not depend on the orderings of edges on cycles. Hence it can be shown that 
dual 
matroids cannot  distinguish triangulations of homology spheres\footnote{These are compact 
connected 3-manifolds whose homology groups are trivial. Unlike in the 2-dimensional case, this does 
not imply that the fundamental group is trivial.} from triangulations of the 3-sphere. 
While the later ones are always embeddable, this is not true for triangulations of homology 
spheres in general. 
Thus we restrict our attention to simply connected simplicial complexes.
Under these necessary restrictions we obtain the following 3-dimensional analogue of Whitney's 
theorem. 

\begin{thm}\label{Whitney_intro}\label{Whitney}
 Let $C$ be a simply connected 2-dimensional  simplicial complex whose dual matroid $M$ is local. 

Then $C$ is embeddable in 3-space if and only if $M$ is graphic.
\end{thm}

Tutte's characterisation of graphic matroids \cite{TutteBook} yields the following 
consequence.

\begin{cor}\label{tutte_cor}
  Let $C$ be a simply connected simplicial complex whose dual matroid $M$ is local.

Then $C$ is embeddable in 3-space if and only if $M$ has no minor isomorphic to 
$U^2_4$, the fano plane, the dual of the fano plane or the duals of either $M(K_5)$ or $M(K_{3,3})$.
\qed
\end{cor}

\vspace{.3 cm}

We further apply dual matroids to study embeddings in 3-space of -- not necessarily simply 
connected -- simplicial complexes with locally small separators as follows.

Given a 2-dimensional simplicial complex $C$, the \emph{link graph}, denoted by $L(v)$, at a 
vertex $v$ of $C$ is 
the graph whose vertices are 
the edges incident with $v$ and whose edges are the faces incident with $v$ and their incidence relation is as in $C$. 
If the link graph at $v$ is not connected, we can split $v$ into one vertex for each connected 
component. There is a similar splitting operation at edges of $C$.
It can be shown that no matter in which order one does all these splittings, one always ends up 
with the same simplicial complex, \emph{the split complex of $C$}.

\begin{figure} [htpb]   
\begin{center}
   	  \includegraphics[height=3cm]{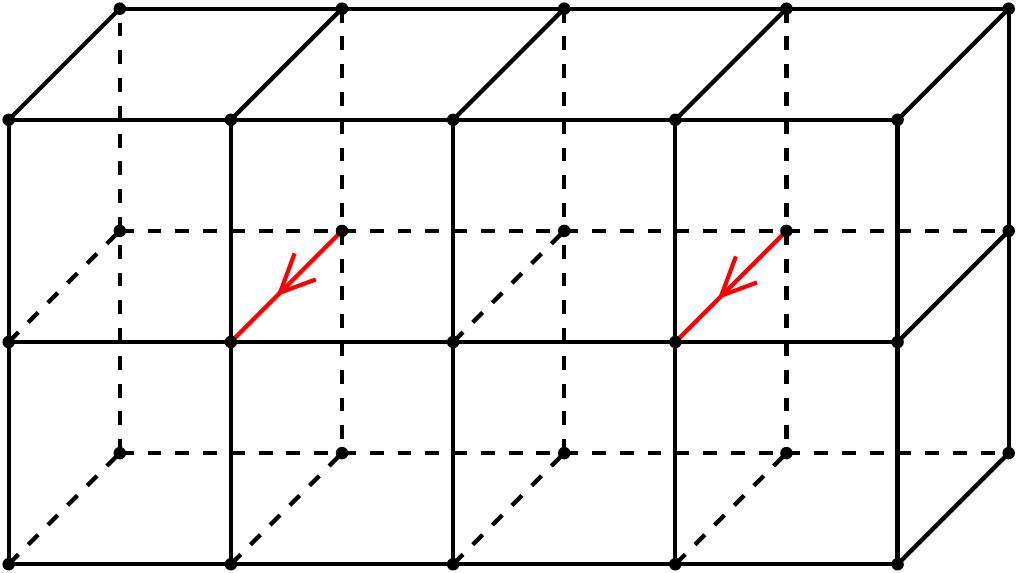}
   	  \caption{The $4\times 2 \times 1$-grid whose faces are the 4-cycles. It can be shown 
that the complex obtained by identifying the two edges coloured red cannot be embedded in 3-space. 
}\label{fig:3grid}
\end{center}
   \end{figure}

It can be shown that if a simplicial complex embeds  topologically into \Sthree, then so does its 
split complexes.
However, the converse is not true. For an example see \autoref{fig:3grid}. 
Here we give a characterisation of when certain simplicial complexes embed, where 
one of the conditions is that the split complex embeds.

\begin{thm}\label{embed_via_matroid}
 Let $C$ be a globally 3-connected simplicial complex and $\hat C$ be its split complex. 
Then $C$ embeds into \Sthree\ if and 
only if $\hat C$ embeds into \Sthree\ 
and the dual matroid of $C$ is the cycle matroid of a graph $G$ and for any vertex or edge of $C$
the set of faces incident with it is a connected edge 
set of $G$. 
\end{thm}
Here a simplicial complex $C$ is \emph{globally 3-connected}\footnote{In \autoref{appendixA} we 
give 
an equivalent definition 
directly in 
terms of $C$.} if its dual matroid is 3-connected. 
For an extension of \autoref{embed_via_matroid} to simplicial complexes that are not globally 
3-connected, see \autoref{embed_via_graph} below. 
\begin{comment}
 Where does the global 3-connectivity come in? The graph G need not be the dual graph a rotation 
system of the split complex by 
assumption. However it will be by uniqueness of the graph. 
\end{comment}

\begin{figure} [htpb]   
\begin{center}
   	  \includegraphics[height=1.5cm]{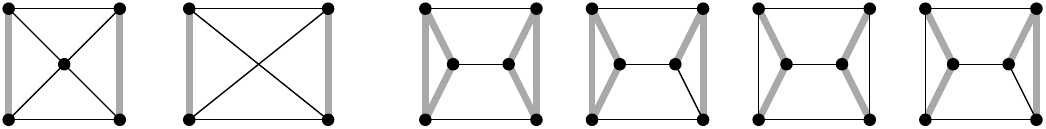}
   	  \caption{The six obstructions characterising connectedness of $X$. In these graphs we 
depicted the edge set $X$ in grey.}\label{fig:constraint}
\end{center}
   \end{figure}

The condition that a given set of elements of the dual matroid is connected (in some graph 
representing that matroid) can be 
characterised by a 
finite list of obstructions as follows. 
Given a matroid $M$ and a set $X$ of its elements, a \emph{constraint minor} of $(M,X)$ is obtained 
by contracting arbitrary elements or deleting elements not in $X$. 
In \cite{3space3}, we prove for any 3-connected graphic matroid $M$ (that is 
a 3-connected graph) with an edge set $X$ that $X$ is connected in $M$ if and only if $(M,X)$ has 
no constraint minor from the 
finite list depicted in \autoref{fig:constraint}. 

\vspace{.3 cm}

In \cite{3space1}, we introduced \emph{space minors} of simplicial complexes and proved 
that a simply connected locally 3-connected simplicial complex $C$ embeds in 3-space if and only if 
it does not have a 
space minor from a finite list $\Lcal$ of obstructions. Using \autoref{embed_via_matroid} we 
can further extend this characterisation from simply connected simplicial complexes to those whose 
split complex is simply connected. 

\begin{thm}\label{combined}
  Let $C$ be a globally 3-connected simplicial complex such that the split complex is 
simply connected and locally 3-connected\footnote{In \cite{3space5} we discuss how this 
result can be extended to simplicial complexes whose split complexes are not local 3-connected.}. 
Then $C$ embeds into \Sthree\ if and 
only if its split complex has no space minor from $\Lcal$ and the 
dual matroid has no constraint minor from the list of \autoref{fig:constraint}. 
\end{thm}

If we do not require global 3-connectivity in \autoref{combined}, there are infinitely many 
obstructions to embeddability, see \autoref{sec:examples}.
We remark that \autoref{Whitney_intro} can be extended from simply connected simplicial 
complexes to those whose split complex is simply connected. 

\vspace{.3 cm}

The paper is structured as follows. In \autoref{sec:dual_matroids} we prove 
\autoref{3dual_matroid}, which is used in the proof of \autoref{Whitney_intro} and 
\autoref{embed_via_matroid}. In \autoref{secW} we prove \autoref{Whitney_intro}. In \autoref{sec4} 
we prove \autoref{embed_via_matroid} and \autoref{combined}. Finally in \autoref{sec:examples} we 
construct infinitely 
many obstructions to embeddability in 3-space (inside the class of simplicial complexes with a 
simply connected and locally 3-connected split 
complex). 

\vspace{.3 cm}

For graph we follow the notations of \cite{DiestelBookCurrent} and for matroids \cite{oxley2}. 
Beyond 
that we rely on some definitions of \cite{3space2}.

\section{Dual matroids}\label{sec:dual_matroids}

In this section we prove \autoref{3dual_matroid} and the fact that a simplicial complex and its 
split complexes have the same dual matroid, 
which are used in the proofs of \autoref{Whitney_intro} and \autoref{embed_via_matroid}. 

A \emph{directed simplicial complex } is a simplicial complex $C$ together with an assignment of a 
direction to each 
edge of $C$ and together with an assignment of a cyclic orientation to each face of $C$. A 
\emph{signed incidence vector} of an edge $e$ of $C$ has one entry for every face $f$; this entry 
is zero if $e$ is not incident with $f$, it is plus one if $f$ traverses $e$ positively and minus 
one otherwise.  

The matrix given by all signed incidence vectors is called the \emph{(signed) edge/face incidence 
matrix}.
The \emph{dual 
matroid} of a simplicial complex is the matroid represented by the edge/face incidence matrix of 
$C$ 
over 
the finite field $\Fbb_3$. 

Although in this paper we work with directed simplicial complexes, dual matroids do not depend on 
the 
chosen directions. Indeed,  changing a direction of an edge or of a face of $C$ changes the linear 
 representation  of the dual matroid but not the matroid itself.

A matrix $A$ is a \emph{regular representation} (or \emph{representation over the integers}) of a 
matroid $M$ if all its entries are integers and the columns are indexed with the 
elements of $M$.
Furthermore for every circuit $o$ of $M$ there is a $\{0,-1,+1\}$-valued vector\footnote{A 
\emph{vector} is an element of a vector space 
$k^S$, where $k$ is a field and $S$ is a set. In a slight abuse of notation, in this paper we also 
call elements of modules of the form 
$\Zbb^S$ vectors.} $v_o$ in the span over $\Zbb$ of the rows of $A$ whose support is $o$. 
And the vectors $v_o$ span over $\Zbb$ all row vectors of $A$.

\subsection{Proof of \autoref{3dual_matroid}}

Let $C$ be a directed simplicial complex embedded into $\Sbb^3$, the \emph{dual digraph} of the 
embedding is 
the following. 
Its vertex set is the set of components of $\Sbb^3\sm C$. It has one edge for every face of $C$. 
This face touches one or two components of 
$\Sbb^3\sm C$.
If it touches two components, the edge for that face joins the vertices for these two components. 
The edge is directed from
the vertex whose complement touches the chosen orientation of the face to the other component.
If the face touches just one component, its edge is a loop attached at the vertex corresponding to 
that component.
\begin{comment}
 
Note that any two topological embeddings that induce the same planar rotation system 
have the same 
dual digraph. 
\end{comment}

Let $(\sigma(e)|e\in E(C))$ be the planar rotation system of $C$ induced by the topological 
embedding of $C$. 
It is not hard to check that $\sigma(e)$ is a closed trail\footnote{A \emph{trail} is sequence 
$(e_i|i\leq n)$ of distinct edges such that the endvertex of $e_i$ is the starting vertex of 
$e_{i+1}$ for all $i<n$.  A trail is \emph{closed} if the starting vertex of $e_1$ is equal to 
the endvertex of $e_n$.}  in the dual graph. 
The \emph{dual complex} of the embedding is the directed simplicial complex obtained from the dual 
digraph 
by adding for each edge of $C$ the 
cyclic orderings of the cyclic orientations $\sigma(e)$ as faces and we choose their orientations 
to 
be $\sigma(e)$. 

\begin{obs}\label{geo_is_abstract}
 Let $C$ be a connected and locally connected\footnote{A simplicial complex $C$ is \emph{locally 
connected} if all its link graphs are connected. } simplicial complex embedded in \Sthree\ with 
induced 
planar rotation system $\Sigma$. Then the dual complex of the embedding is equal to the dual 
complex 
of $(C,\Sigma)$.
\end{obs}

\begin{proof}
 By {\cite[Lemma 3.4]{{3space2}}}, the local surfaces for 
$(C,\Sigma)$ agree with the local surfaces of the 
embedding\footnote{Local surfaces of embeddings are defined in \cite{3space2}.}. Hence these two 
complexes have the vertex set. As they also have the same incidence 
relations between edges and vertices and edges and faces, they must coincide.
\end{proof}

By \autoref{geo_is_abstract} and the definition of `generated over the integers' and by 
\autoref{char_reg_represent}, in order to prove \autoref{3dual_matroid} it suffices to 
show that the dual complex for  $(C,\Sigma)$ is nullhomologous\footnote{A simplicial complex $C$ 
is \emph{nullhomologous} if the face boundaries of $C$ generate all cycles over the integers. This 
is equivalent to the condition that the face boundaries of $C$ generate all cycles over the field 
$\Fbb_p$ for every prime $p$.}. 

First we prove this in the special case when $C$ is nullhomologous and locally 
connected. 

\begin{lem}\label{dual_nullh}
 Let $C$ be a nullhomologous locally connected simplicial complex together with a planar rotation 
system 
$\Sigma$ such that local surfaces for $(C,\Sigma)$ are spheres.\footnote{This last 
property follows from the first two if we additionally assume that $\Sigma$ is induced by a 
topological embedding in \Sthree\ by {\cite[Theorem 6.1]{{3space2}}}.} Then 
the dual complex $D$ 
of $(C,\Sigma)$ is nullhomologous. 
\end{lem}

\begin{proof}
 By {\cite[Lemmas 6.3,6.5,6.7]{{3space2}}} the complexes $C$ and 
$D$ satisfy euler's formula, that is:
 \[
  |V(C)|-|E|+|F|-|V(D)|=0
 \]
Hence we deduce that $D$ nullhomologous by applying the `Moreover'-part of 
{\cite[Lemma 6.3]{{3space2}}} for every prime $p$. 
\end{proof}

Next we shall extend  \autoref{dual_nullh} to simplicial complexes that are only locally connected. 

\begin{lem}\label{dual_nullh_stronger}
 Let $C$ be a locally connected simplicial complex together with a planar rotation 
system $\Sigma$ that is induced by a topological embedding $\iota$ in \Sthree.  Then the dual 
complex $D$ 
of $(C,\Sigma)$ is nullhomologous. 
\end{lem}

\begin{proof}
 By {\cite[Theorem 7.1]{{3space2}}}
there is a simplicial complex $C'$ that is obtained from $C$ by subdividing edges, baricentric 
subdivisions of faces and adding faces 
along closed trails. And $C'$ is nullhomotopic and has an embedding $\iota'$ 
into \Sthree\ that induces\footnote{This means that we obtain $\iota$ from $\iota'$ by deleting 
the newly added faces, contracting the newly added subdivision edges and undoing the baricentric 
subdivisions.} $\iota$. 
Let $D'$ be the dual of $\iota'$. 
By \autoref{dual_nullh}, $D'$ is nullhomologous. 

We shall deduce that $D$ is nullhomologous by showing that reversing each of the operations in the 
construction of $C'$ from $C$ preserves being nullhomologous in the dual. We call such an 
operation \emph{preserving}.

\begin{sublem}\label{l1}
Subdividing an edge is preserving.
\end{sublem}

\begin{proof}
Subdividing an edge in the primal corresponds to adding a copy of a face in the dual. Clearly, the 
deletion of the copy preserves being nullhomologous for the dual. 
\end{proof}

\begin{sublem}\label{l2}
A baricentric subdivision of a face is preserving. 
\end{sublem}

\begin{proof}
It suffices to show that the subdivision by a single edge is preserving.
Subdividing a face by an edge in the primal corresponds to replacing an edge in the dual by two 
edges in parallel and adding a face containing precisely these two edges. 
Reversing this operation preserves being  nullhomologous. 
\end{proof}

\begin{sublem}\label{l3}
Adding a face is preserving. 
\end{sublem}

\begin{proof}
Adding a face in the primal corresponds to coadding\footnote{A complex $A$ is obtained 
from a complex $A'$ by \emph{coadding} an edge $e$ if $A'$ is obtained from $A$ by contracting 
the edge $e$.} an edge in the dual.
Contracting that edge preserves being  nullhomologous. 
\end{proof}

By \autoref{l1}, \autoref{l2} and \autoref{l3}, the fact that $D'$ is nullhomologous implies that 
$D$ is nullhomologous. 
\end{proof}
 
It remains to prove \autoref{3dual_matroid} for simplicial complexes $C$ that are not 
locally connected. First we need some preparation.

Given a simplicial complex $C$, its \emph{vertical split complex} is obtained from $C$ by replacing 
each 
vertex 
$v$ by one vertex for each 
connected component of $L(v)$, where the edges and faces incident with that vertex are those in its 
connected component.  We refer to these new vertices as the \emph{clones} of $v$. 
 
\begin{obs}\label{is_loc_conS4}
 The vertical split complex of any simplicial complex is locally connected. 
 \qed
\end{obs}

\begin{obs}\label{same_matroid}
 A simplicial complex and its vertical split complex have the same dual matroid.
\end{obs}

\begin{proof}
 A simplicial complex and its vertical split complex have the same edge/face incidence matrix. 
\end{proof}

Given an embedding $\iota$ of a simplicial complex $C$ into \Sthree, we will define what an 
\emph{induced} embedding of the vertical split 
complex is. 

\begin{figure} [htpb]   
\begin{center}
   	  \includegraphics[height=3cm]{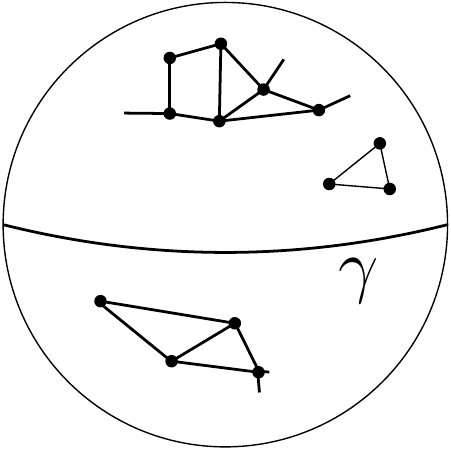}
   	  \caption{The link graph at $v$ embedded into $\partial B$. 
}\label{split_vertex}
\end{center}
   \end{figure}

For that we need some preparation. Let $v$ be a vertex of $C$ whose link graph is not connected.
By changing $\iota$ a little bit locally (but not its induced planar rotation system) if necessary, 
we may assume that there is a 2-ball $B$ 
of small radius around $v$ such that firstly $v$ is the only vertex of $C$ contained in the inside 
of $B$. And secondly its boundary $\partial B$ 
intersects each edge incident with $v$ in a point and each face incident with $v$ in a line. In 
other words, the intersection of $C$ with 
the boundary is the link graph at $v$. As the link graph is disconnected, there is a circle 
(homeomorphic image 
of $\Sbb^1$) $\gamma$ in the 
boundary such that the two components of $B\sm \gamma$ both contain 
vertices of the 
link graph, see \autoref{split_vertex}.

The simplicial complex $C_\gamma$ is obtained from $C$ by replacing the vertex $v$ by two vertices, 
one for each connected component of 
$\partial B\sm \gamma$ that is incident with the edges and faces whose vertices and edges, 
respectively, are in that connected component.

The embedding $\iota$ \emph{induces}\footnote{ The construction of $\iota_C$ depends on the choice 
of $B$. Still we use the term `induced' 
in this context since in this paper we 
consider topological embeddings equivalent if they have the same planar rotation system.  } the 
following embedding $\iota_\gamma$ of 
$C_\gamma$ into \Sthree. 
We pick a disc contained in $B$ with boundary $\gamma$ that intersects $C$ only in $v$. 
We replace $v$ by its two clones -- both with tiny distance from $v$ and one above that disc and 
the other below. We only need to change faces and edges incident with $v$ in a tiny neighbourhood 
around $v$. Faces and edges above and below do not interfere. 

It is easy to see that $\iota$ and $\iota_C$ have the same planar rotation system and that $C$ and 
$C_\gamma$ have the same vertical 
split complex.

A topological embedding of the vertical split complex of $C$ into \Sthree\ is \emph{(vertically) 
induced} by $\iota$ if it is obtained by 
applying the above procedure iteratively until $C_\gamma$ is equal to the vertical split complex of 
$C$. It is 
clear that if $\iota$ is a topological 
embedding of a simplicial complex $C$ into \Sthree, then its vertical split complex has a 
topological embedding into \Sthree\ that is 
induced by $\iota$.

\begin{obs}\label{same_complex}
Let $\iota$ be an embedding of a simplicial complex into \Sthree\ and let $\iota'$ be an induced 
embedding of 
$\iota$ of the vertical split complex.
Then $\iota$ and $\iota'$ have the same dual complex.
\end{obs}

\begin{proof}
In both embeddings, the incidence relation between the local surfaces and the faces is the same. 
Hence both dual complexes have the same 
vertex/edge incidence relation. They also have the same sets of faces as $\iota$ and $\iota'$ have 
the same rotation system. 
\end{proof}

A set $S$ of vertices in a simplicial complex $C$ is a \emph{vertex separator} if $C$ can be 
obtained from two 
disjoint simplicial complexes that each have at 
least one face by gluing them together at the vertex set $S$. As the empty set might also be a 
vertex separator, any simplicial complex 
with no vertex separator is connected. 

\begin{lem}\label{dual_complex}
Let $C$ be a simplicial complex without a vertex separator. Assume that $C$ has an 
embedding $\iota$ into \Sthree. 
Then the dual complex $D$ of $\iota$ is 
nullhomologous. 
\end{lem}

\begin{proof}
 Let $C'$ be the vertical split complex of $C$. 
  By \autoref{is_loc_conS4}, $C'$ is locally connected. 
By assumption $C$ has no vertex separator. Thus $C'$ is connected. 
 Let $\iota'$ be the embedding of $C'$ induced by 
$\iota$. Let $\Sigma'$ be the planar rotation system induced by $\iota'$.

By \autoref{dual_nullh_stronger}, the dual $D'$ for $(C',\Sigma')$ is nullhomologous. By 
\autoref{geo_is_abstract}, $D'$ is the dual complex of $\iota'$. 
By \autoref{same_complex}, $D'$ is equal to $D$. So $D$ is nullhomologous. 
\end{proof}

\begin{lem}\label{conny2}
 Let $C$ be a simplicial complex embedded into \Sthree\ that is obtained from two 
simplicial complexes $C_1$ 
and $C_2$ by 
gluing them together at a set of 
vertices. Assume that $C_2$ has no separating vertex set. 
 Let $G_i$ be the dual graph of the embedding restricted to $C_i$ for $i=1,2$. 
 Then the dual graph of the embedding of $C$ is equal to a graph obtained by gluing together 
$G_1$ 
and $G_2$ at a single vertex.
\end{lem}

\begin{proof}
 We denote the embedding of $C$ into \Sthree\ by $\iota$ and the restricted embedding of 
$C_1$ by $\iota_1$. 
 Suppose for a contradiction that $\iota$ maps interior points of faces of $C_2$ to interior points 
of different local surfaces of 
$\iota_1$.
 Let $\ell$ be a local surface of $\iota_1$ to which an interior point of a face of $C_2$ is mapped 
by $\iota$. 
 Let $C_2'$ be the subcomplex of $C_2$ that contains all faces whose interior points are 
mapped to interior points of $\ell$. Its edges and vertices are those of $C_2$ that are incident 
with these faces. Note that if one interior point of a face is mapped to $\ell$, 
then 
all are. Hence the subcomplex $C_2''$ 
that contains all other faces and their incident vertices and edges contains a face. The 
subcomplexes $C_2'$ and $C_2''$ of $C_2$ can only intersect in points of $C_1$. Hence they only can 
intersect in vertices. Thus $C_2'$ 
and $C_2''$ witness that $C_2$ has a 
separating vertex set contrary to our assumption. 

Thus there is a single local surface of $\iota_1$ to which all interior points of faces of $C_2$ 
are mapped 
by $\iota$. Hence the 
dual graph of $\iota$ is equal to the graph obtained by gluing together $G_1$ and $G_2$ at that 
vertex.
\end{proof}

       \begin{proof}[Proof of \autoref{3dual_matroid}.]
 By applying \autoref{conny2} recursively, we may assume that $C$ has no separating vertex set. 
        Recall that the dual graph of the embedding is the 1-skeleton of the dual 
complex of the embedding. 
       By \autoref{dual_complex}, the edge/face incidence matrix is a representation over the 
integers of the cycle matroid of the dual graph of the embedding. 
       \end{proof}
 
\subsection{Split complexes}
A naive way to define splittings of edges might be to consider the 
incidences at one of their endvertices and split according to that. We shall show that when 
using 
this notion of splitting, split complexes will not have all nice properties we want them to have, 
see \autoref{appendixA}. A more refined definition takes into 
account the incidences at both endvertices, defined as follows.

Given a simplicial complex $C$ and an edge $e$ with two endvertices $v$ and $w$, two faces 
incident with $e$ are \emph{$v$-related} if - when considered as edges of $e$, they have endvertices 
in the same connected component of the link graph $L(v)-e$ with the vertex $e$ removed. 
Analogously, we define \emph{$w$-related}. 
Two faces $f_1$ and $f_2$ incident with $e$ are in the \emph{same connected 
component at $e$} if there is a chain of faces incident with $e$ from $f_1$ to $f_2$ such that 
adjacent faces in the chain are $v$-related or $w$-related. Note that `being in the same connected 
component at $e$' is the equivalence relation generated from the union of `$v$-related' and 
`$w$'-related.

The simplicial complex obtained from $C$ by \emph{splitting} the edge $e$ is obtained from 
$C$ by replacing the edge $e$ by 
one copy $e_X$ for every connected component $X$ at $e$. The faces incident with $e_X$ are those in 
$X$.

We refer to the edges $e_X$ as the \emph{clones} of $e$. If we apply several splittings, we extend 
the 
notion of cloning iteratively so that each edge of the resulting simplicial complex is cloned from 
a 
unique edge of $C$.

If we split an edge in a nontrivial way, then the resulting simplicial complex has the same number 
of faces 
but at least 
one edge more. As in a simplicial 
complex every 
edge is incident with a face, we can only split edges a bounded number of times. A simplicial 
complex obtained from $C$ by splitting 
edges 
such that for every edge there is only one component at $e$ is called an \emph{edge split complex} 
of $C$. As explained above, every simplicial complex has an edge split complex. 

Since splitting edges, does not change the 2-blocks of the link graphs, splittings of edges commute.
In particular, edge split complexes are unique. In 
the following we will talk about 
`the edge split complex'. 

The \emph{split complex} of a simplicial complex $C$ is the vertical split complex of its edge 
split complex. Clearly, splitting a vertex does not change the edge split complex. 

 \begin{eg}\label{egx}
A simplicial complex, its vertical split complex and its edge split complex have the same split 
complex. Locally 
2-connected\footnote{A simplicial complex is \emph{locally 2-connected} if its link graphs are 
connected and have no 
cutvertices.} simplicial complexes are equal to 
their split complex. 
 \end{eg}

\begin{lem}\label{same_matroid2_edge}
 A simplicial complex and its edge split complex have the same dual matroid.
\end{lem}

\begin{proof}
We shall show that a simplicial complex $C$ and a simplicial
complex $C'$ have the same dual matroid, where we obtain $C'$ from $C$ by splitting an edge $e$.
Once this is shown, the lemma follows inductively as an edge split complex is obtained by a 
sequence 
of edge splittings.

Clearly, $C$ and $C'$ have the same set of faces. Hence their dual matroids have the same ground 
sets. 

 The vectors indexed by clones of the edge $e$ of the edge/face incidence matrix 
$A'$ of $C'$ sum up to the vector indexed by $e$ of the edge/face incidence matrix $A$ of $C$.
 Hence the vectors indexed by edges of $A'$ generate the vectors indexed by edges of $A$.
 So it remains to show that any vector indexed by a clone $e'$ of $e$ of $A'$ is 
generated by the vectors indexed by edges of 
$A$. 

Let $v$ be an endvertex of $e$. 
Let $K$ be the connected component of the link graph $L(v)$ of $C$ at $v$ that contains $e$. 
Let $Y$ be the union of the 
components $Y'$ of $K-e$ such that faces incident with $e'$ -- when considered as edges of $L(v)$ 
-- have an endvertex in $Y'$. 
The sum over all vectors indexed by edges $y\in V(Y)$ of $A$ is the vector indexed by $e'$ of 
$A'$.
Since $e'$ was an arbitrary clone, the vectors indexed by edges of $A$ generate the vectors 
indexed by edges of $A'$.

We have shown that splitting a single edge preserves the dual matroid. Since the edge split complex 
is obtained by splitting edges, 
it must 
have the same dual matroid as the original complex. 
\end{proof}

\begin{cor}\label{same_matroid2}
 A simplicial complex and its split complex have the same dual matroid.
\end{cor}

\begin{proof}
A simplicial complex and its vertical split complex have the same incidence relations between edges 
and faces. Hence this is a consequence 
of \autoref{same_matroid2_edge}. 
\end{proof}

\section{A Whitney type theorem}\label{secW}

In this section we prove \autoref{Whitney_intro}. 

In general the dual matroid of a simplicial complex $C$ does not contain enough information to 
decide 
whether $C$ is embeddable in 3-space. For example, the dual matroid of the cone over $K_5$ 
consists 
of a bunch of loops. So it cannot distinguish this non-embeddable simplicial complex from other 
embeddable 
ones. The following fact gives an explanation of this phenomenon (in the notation of that 
fact: from the graph $G$ we can in general not reconstruct the matroid $M[v]$).
Given a vertex $v$ of a simplicial complex, we denote the dual matroid of the link graph at $v$ by 
$M[v]$. 

\begin{fact}\label{not_reconstruct}
Let $C$ be a simplicial complex embedded in $\Sbb^3$. Then the dual matroid $M$ restricted to the 
faces 
incident with $v$ is represented by a graph $G$. Moreover, $G$ can be obtained from some graph 
representing $M[v]$ by identifying vertices.  
\end{fact}

\begin{proof}
By \autoref{3dual_matroid} $M$ is the cycle matroid of the dual graph of the embedding of $C$. So 
$G$ is the restriction of that graph to the faces incident with $v$. 

By $G'$ be denote the `local 
dual graph' of $C$ at $v$. This is defined as the `dual graph' but with `$\Sbb^3$' replaced by `a 
small neighbourhood $U$ around $v$' in the embedding. Clearly, $G'$ represents $M[v]$. 
We obtain the vertices of $G$ 
from those of $G'$ by identifying those vertices for components of $U\sm C$ that lie in the same 
component of $\Sbb^3\sm C$. 
The `Moreover'-part follows. 
\end{proof}

To exclude the phenomenon described in \autoref{not_reconstruct} we restrict our attention to 
simplicial complexes $C$ 
whose dual matroid captures the local structure at all vertices of $C$, defined as follows. 
Given a  simplicial complex $C$ with dual matroid $M$, we say that $M$ is \emph{local} if
for every vertex $v$ the matroid $M[v]$ is equal to $M$ restricted to the faces incident with $v$.

Furthermore matroids (of graphs and also of 
simplicial complexes) do not depend on the orderings of edges on cycles. Hence it can be shown that 
dual 
matroids cannot  distinguish triangulations of homology spheres\footnote{These are compact 
connected 3-manifolds whose homology groups are trivial. Unlike in the 2-dimensional case, this 
does 
not imply that the fundamental group is trivial.} from triangulations of the 3-sphere. 
While the later ones are always embeddable, this is not true for triangulations of homology 
spheres. 
Thus we restrict our attention to simply connected simplicial complexes.

If we exclude these two phenomenons,
\autoref{Whitney}, stated in the Introduction, characterises when a simplicial complex is 
embeddable just in terms of its dual matroid.

\begin{rem}
 The assumptions of \autoref{Whitney} can be interpreted as some face maximality assumption.
By {\cite[Theorem 7.1]{{3space2}}} this is true for being simply connected. 
For locality, let $C$ be any embeddable simplicial complex embeddable. By 
\autoref{not_reconstruct} we can add faces until for every vertex $v$ the matroid $M[v]$ is equal 
to 
$M$ restricted to the faces incident with $v$. This preserves being simply connected.
\end{rem}

Now we prepare for the proof of \autoref{Whitney}. 
\begin{lem}\label{eulerian_structure_exists}
Let $H$ be a graph whose cycle matroid is the dual matroid $M$ of a simplicial complex $C$. 
There is a directed graph $\vec{H}$ with underlying graph $H$ such that for all edges $e$ of $C$ 
the signed vectors are 3-flows\footnote{A \emph{3-flow} in a directed graph $\vec{H}$ is an 
assignment of integers to the edges of $\vec{H}$ that satisfies Kirchhoff's first law modulo three  
at every vertex of $\vec{H}$.}.
\end{lem}

\begin{proof}
First we consider the case when $H$ is 2-connected. 
We start with an arbitrarily directed graph $\vec{H}$ with underlying graph $H$ some of whose 
directions of the edges we might reverse later on in the argument.
Since $H$ is 2-connected, the set of edges incident with a vertex is a bond of $H$, which is called 
the \emph{atomic bond} of $v$.
By elementary properties of representations, there is a vector $b_v$ with all entries 
$-1$, $+1$ or $0$ that has the same support\footnote{The \emph{support} of a vector is the set of 
coordinates with nonzero values. } as the atomic bond at $v$.

Given an edge $e$ of $H$ and one of its endvertices $v$, we say that $e$ is 
\emph{effectively directed towards} $v$ with respect to a vector $b$ with entries in $\Zbb$ if 
$\vec{e}$ is directed towards $v$ and $b(e)$ is positive or $\vec{e}$ is directed away from 
$v$ and $b(e)$ is negative. 
First we shall prove that we can modify the directions of the edges of $\vec{H}$ such that all 
edges 
$e$ of $H$ are directed such that for some endvertex $v$ they are effectively directed 
towards $v$ with respect to the at $b_v$.

Let $T$ be a spanning tree of $H$. 
Since $T$ does not contain any cycle, we can pick the $b_v$ such that if $vw$ is an edge of $T$, 
then $b_v(vw)=-b_w(vw)$. 
Hence an edge $vw$ of $T$ is effectively directed towards $v$ with respect to $b_v$ if and only if 
it is effectively directed towards $w$ with respect to $b_w$. 
So by reversing the direction of an edge if necessary\footnote{To be very formal, we delete the 
edge 
from the graph and glue it back the other way round. Note that we do not change the director.}, we 
may assume that every edge $vw$ of $T$ is effectively directed towards $v$ with respect to $b_v$ 
and 
also
effectively directed towards $w$ with respect to $b_w$. 

Next let $xy$ be an edge not in $T$. By reversing the direction of $xy$ if necessary we may assume 
that $xy$ is effectively directed towards $x$ with respect to $b_x$. Our aim is to show that $xy$ 
is 
effectively directed towards $y$ with respect to $b_y$. 
Let $C$ be the fundamental circuit of $xy$ with respect to $T$. 
By elementary properties of representations, there is a vector $v_C$ with support $C$ that is 
orthogonal over 
$\Fbb_3$ to all the vectors $b_z$ for vertices $z$ on $C$. At all vertices $z$ of $C$ 
except possibly $y$, the two edges on $C$ incident with $z$ are 
effectively directed towards $z$ with respect to the vector $b_z$. 
Hence for $v_C$ to be orthogonal, precisely one of these edges must be effectively directed towards 
$z$ with respect to $v_C$.
Using this property inductively along $C$, we deduce that of the two edges on $C$ incident with $y$ 
also precisely one is effectively directed towards $y$ with respect to $v_C$.
Since $b_y$ is orthogonal to $v_C$ and the edge incident with $y$ that is on $T$ and $C$ is 
effectively directed towards $y$ with respect to $b_y$, also $xy$
must be effectively directed towards $y$ with respect to $b_y$.

Hence our final directed graph $\vec{H}$ has the property that all edges $e$ of $H$ are effectively 
directed towards any of their endvertices $v$ with respect to $b_v$. 
Since signed vectors of edges $e$ of $C$ are orthogonal at to $b_v$, it follows 
that it accumulates 0 (mod 3) at all vertices $v$. So the signed vectors of $C$ are 3-flows for 
$\vec{H}$. 
This completes the proof if $H$ is 2-connected. 
If $H$ is not 2-connected, we do the same construction independently in every 2-connected component 
and the result follows. 
\end{proof}

First we prove \autoref{Whitney} under the additional assumption that $C$ is locally 2-connected: 

\begin{lem}\label{Whitney2}
 Let $C$ be a simply connected locally 2-connected simplicial complex whose dual matroid is local. 

Then $C$ is embeddable in 3-space if and only if $M$ is graphic.
\end{lem}

\begin{proof}
Assume that $C$ is embeddable and let $D$ by its dual complex. Then by \autoref{3dual_matroid} $M$ 
is equal to the cycle matroid of the 1-skeleton of $D$. In particular $M$ is graphic.

Now conversely assume that $C$ is a simply connected simplicial complex such for every vertex $v$ 
the 
matroid $M[v]$ is equal to dual matroid $M$ restricted to the faces incident with $v$; and that 
there is a graph $G$ whose cycle matroid is $M$. We pick an 
arbitrary direction at each edge of $C$ and an arbitrary orientation at each face of $C$. 
Our aim is to construct a planar rotation system $\Sigma$ of $C$ and apply 
{\cite[Theorem 1.1]{{3space2}}} to deduce that 
$C$ is embeddable.  
 
By \autoref{eulerian_structure_exists} there is a direction $\vec{G}$ of $G$ such that the signed 
incidence vector $v_e$ for each edge $e$ of $C$ is a 3-flow in  $\vec{G}$. As the link 
graph $L(v)$ at each 
vertex $v$ is 2-connected, none of its vertices $e$ is a cutvertex. Hence the edges incident with 
$e$ in $L(v)$ form a bond. So they form a circuit in the dual matroid $M[v]$. Thus by assumption 
the support of $v_e$ is a circuit in the matroid $M$. By the construction of $\vec{G}$, the signed 
vector $v_e$ is a directed cycle\footnote{A vector $v$ whose entries are in $\{0,+1,-1\}$ is a 
\emph{directed cycle} if its support is a cycle and it satisfies Kirchhoff's first law at every 
vertex, see \cite{DiestelBookCurrent}.} in 
$\vec{G}$. 
This directed cycle defines a cyclic orientation $\sigma(e)$. In terms of $C$ this is a cyclic 
orientation of the oriented faces incident with the directed edge $\vec{e}$. 
Put another way $\Sigma=(\sigma(e)|e\in E(C))$ is a rotation system.

Our aim is to prove that $\Sigma$ is planar. So let $v$ be a vertex of $C$ and let $\Sigma_v$ be 
the rotation system of the link graph $L(v)$ induced by $\Sigma$. This rotation system of $L(v)$ 
defines an embedding of $L(v)$ in a 2-dimensional oriented surface 
$S_v$ in the sense of \cite{moharThomassen}\footnote{This is explained in more detail in 
\cite{3space2}. }. 
It remains to show the following.

\begin{sublem}\label{euler_genus}
$S_v$ is a sphere. 
\end{sublem}

\begin{proof}
As the graph $L=L(v)$ is connected, $S_v$ is connected. Thus it suffices to show that it has Euler 
genus two, that is:
\begin{equation}\label{euler}
 V_L-E_L+F_L=2
\end{equation}
Here we abbreviate: $|V(L)|=V_L$, $|E(L)|=E_L$ and $F_L$ denotes the faces of the embedding of 
$L(v)$ in $S_v$. 

We denote the dual graph of the embedding of $L$ in $S_v$ by $H$.
Our aim is to show that $H$ is equal to the restriction $R$ of $G$ to the faces incident with $v$.
We obtain $S'$ from $R$ by gluing on each directed cycle $v_e$ the face 
$\sigma(e)$. Similarly as in \cite{3space2} we use the Edmonds-Hefter-Ringel 
rotation principle {\cite[Theorem 3.2.4]{{moharThomassen}}} to deduce that $L$ is the 
surface dual of $R$ with respect to the embedding into $S'$.
In particular $S'=S$ and $R$ is equal to $H$.

Having shown that $R$ is the surface dual of $L$, we conclude our proof of \autoref{euler} as 
follows. 
We denote the dimension of the cycle space of $L$ by $d$. We have $V_L-E_L=-d+1$ and $F_L=V_R$ 
(where $V_R$ is the number of vertices of $V_R$). Hence in order to prove \autoref{euler} it 
suffices to show that $d=V_R-1$. This follows from the assumption that the cycle 
matroid of $R$ is the dual of the cycle matroid of $L$. Indeed, the cycle 
matroid of $L$ is 2-connected by assumption. 
\end{proof}
\end{proof}

\begin{proof}[Proof of \autoref{Whitney}.] 
As in the proof of \autoref{Whitney2}, by \autoref{3dual_matroid} 
it suffices to show that any simply connected simplicial complex $C$ whose dual matroid $M$ is 
graphic and local can be embedded in 3-space.

We prove this in two steps. First we prove it for locally connected simplicial complexes. 
We prove this by induction. The base case is when $C$ is locally 2-connected and this is dealt with 
in  \autoref{Whitney2}. So now we assume that $C$ has a vertex $v$ such that the link graph $L(v)$ 
has a cut vertex\footnote{A vertex $v$ of a graph is a \emph{cut vertex} if the component 
of the graph containing $v$ with 
$v$ removed is disconnected.}; and that we proved the statement for every simplicial complex as 
above such that 
it has a 
fewer number of cutvertices -- summed over all link graphs. Let $e$ be an edge of $C$ that is a 
cutvertex in $L(v)$. 

\begin{sublem}\label{identi}
The simplicial complex $C$ is obtained from a simplicial complex $C'$ by identifying two 
vertex-disjoint edges $e_1$ and $e_2$ onto $e$. 
\end{sublem}

\begin{proof} In the link graph $L(v)$, let $f_1$ and $f_2$ be two edges 
incident with $e$ that are in different 2-blocks of $L(v)$. 
Hence $L(v)$ has a 1-separation $(X_1,X_2)$ with cutvertex $e$ such that $f_i$ is in the 
side $X_i$ for $i=1,2$.

 Let $w$ be the endvertex of $e$ in $C$ different from $v$. Our aim is to construct a 
 1-separation $(Y_1,Y_2)$ with cutvertex $e$ of $L(w)$ such that $X_i$ and $Y_i$ agree when 
restricted to the edges incident with $e$ for $i=1,2$. For that we have to show that if two such 
edges are in different $X_i$ then they do not lie in the same 2-block of $L(w)$. That is, in the 
matroid $M[w]$ they do not lie in a common circuit consisting of edges incident with $e$. By the 
assumption, this property is true in $M[w]$ if and only if it is true in $M$ if and only if it is 
true in $M[v]$, which it is not true as  $(X_1,X_2)$ is a 1-separation.

We obtain $C'$ from $C$ by replacing $v$ by two new vertices $v_1$ and $v_2$ and $w$ by two new 
vertices $w_1$ and $w_2$. A face or edge incident with $v$ is in $v_i$ if and only if it is in 
$X_i$. Similarly, a face or edge incident with $w$ is incident with $w_i$ if and only if it is in 
$Y_i$. Thus every edge or face incident with $v$ is incident with precisely one of $v_1$ and $v_2$ 
except for the edge $e$ for which we introduce two copies, which we denote by $e_1$ and $e_2$.
The same is holds with `$w$' in place of `$v$'. Clearly, the edge $e_i$ joins $v_i$ and 
$w_i$.  Hence 
$C'$ has the desired properties. 
\end{proof}

\begin{sublem}\label{X}
The edges $e_1$ and $e_2$ lie in different connected components of $C'$. 
\end{sublem}

\begin{proof}
The simplicial complex  $C/e$ is simply connected and obtained from $C'/\{e_1,e_2\}$ by identifying 
the vertices $e_1$ and $e_2$ onto $e$. 
Since $C/e$ is not locally connected at $e$ we can apply 
{\cite[Lemma 5.1]{{3space2}}} to deduce that $e$ has to be a cutvertex of $C/e$. 

Since the link graph $L(e)$ of $C/e$ is a disjoint union of the connected link graphs $L(e_1)$ and 
$L(e_2)$ of $C'/\{e_1,e_2\}$, two faces incident with the same edge $e_i$ in $C'$ cannot be cut off 
by $e$ in $C/e$. Hence the only way $e$ can cut $C/e$ is that $e_1$ and $e_2$ are cut off from one 
another. Put another way, $e_1$ and $e_2$ lie in different connected components of $C'$. 
\end{proof}

For $i=1,2$, let $C_i$ be the component of $C'$ containing $e_i$ and $M_i$ the dual matroid of 
$C_i$. We may assume that $C$ is connected. Hence $C'$ is the disjoint union of the $C_i$.
By \autoref{X} and \autoref{same_matroid2_edge}, the dual matroid $M$ of $S$ is the disjoint 
union of the matroids $M_i$. So we can 
apply the induction hypothesis to each simplicial complex $C_i$. So all $C_i$ are embeddable.  
Analogously to {\cite[Lemma 5.2]{{3space2}}} one proves that $C$ is embeddable in 
3-space\footnote{An alternative is the following: it is 
easy to see that a 
simplicial complex $S$ is 
embeddable if and only if $S/e$ is embeddable for some nonloop $e$. So the $C_i/e_i$ are 
embeddable. Then 
by {\cite[Lemma 5.2]{{3space2}}} $C/e$ is 
embeddable. So $C$ is embeddable.}.

Finally, we prove the statement for arbitrary simplicial complexes.
Again, we prove it by induction. This time the locally connected case is the base case. 
So now we assume that $C$ has a vertex $v$ such that the link graph $L(v)$ is disconnected; and 
that we proved the statement for every simplicial complex as above such that the number of 
components of link 
graphs minus the total number of link graphs is smaller. As $C$ is simply connected, by 
{\cite[Lemma 5.1]{{3space2}}} the vertex $v$ is a 
cutvertex of $C$. 
That is, $C$ is obtained from gluing together two simplicial complexes $C'$ and $C''$ at the vertex 
$v$. Since splitting vertices preserves dual matroids, the dual matroid of $C$ is 
the disjoint union of the dual matroid of $C'$ 
and the dual matroid of $C''$. Thus the simplicial complexes $C'$ and $C''$ are embeddable in 
$\Sbb^3$ by induction. Hence by 
{\cite[Lemma 5.2]{{3space2}}} $C$ is embeddable.
\end{proof}

\begin{rem}
The proof of \autoref{Whitney} works also if we change the definition of dual matroid in that we 
replace `$\Fbb_3$' by `$\Fbb_p$ with $p$ prime and $p>2$'. By \autoref{3dual_matroid}, if $C$ is 
embeddable, 
the signed incidence vectors of the edges of $C$ 
generate the same matroid over any field $\Fbb_p$ with $p$ prime. So if $C$ is embeddable all these 
definitions of dual matroids coincide. 

The special role of $p=2$ is visible in \autoref{tutte_cor}, where we have to exclude the 
matroid $U_{2,4}$, which is representable over any field $\Fbb_p$ with $p$ prime and $p>2$ but not 
over $\Fbb_2$.
\end{rem}

\section{Constructing embeddings from embeddings of split complexes}\label{sec4}

In this section we prove \autoref{embed_via_matroid}. We subdivide this proof in four subsections. 

\subsection{Constructing embeddings from vertical split complexes}

\begin{lem}\label{loc_surface_identify}
 Let $C$ be a simplicial complex obtained from a simplicial complex $C'$ by identifying two vertices 
$v$ and $w$.
 Let $\iota'$ be a topological embedding of $C'$ into \Sthree. Assume that there is a local surface 
of $\iota'$ that contains both $v$ and 
$w$. Then there is a topological embedding of $C$ into \Sthree\ that has the same dual graph as 
$\iota'$.
\end{lem}

\begin{proof}
We join $v$ and $w$ by a copy of the unit interval $I$ inside the local surface of $\iota'$ that 
contains them both. We may assume that there is an open cylinder around $I$ that does not 
intersect $C'$. 
We obtain a topological embedding $\iota$ of $C$ from $\iota'$ by moving $v$ 
along $I$ to $w$. We do this in such a way that we change the edges and faces incident with 
$v$ only inside the small cylinder. It is clear that $\iota'$ and $\iota$ have the same dual 
graph. 
\end{proof}

\begin{lem}\label{is_conny}
Let $x$ be a vertex or edge of a simplicial complex $C$ embedded into \Sthree. The set of faces 
incident with 
$x$ is a connected edge set of the 
dual graph 
of the embedding.
\end{lem}

\begin{proof}
If $x$ is an edge, then the set of faces incident with $x$ is a closed trail, and hence connected.
Hence it remains to consider the case that $x$ is a vertex. 
Let $H_x$ be the dual graph of the link graph at $x$ with respect to the embedding in the 2-sphere 
given by the embedding of $C$. The restriction $R_x$ of the dual graph of the embedding of $C$ to 
the faces 
incident with $x$ is obtained from $H_x$ by identifying vertices. Since $H_x$ is connected, also 
$R_x$ is connected. This completes the proof.
\end{proof}

Given a simplicial complex $C$ and a topological embedding $\iota$ of its vertical split complex 
into \Sthree, 
we say that $\iota$ satisfies the 
\emph{vertical dual graph connectivity constraints} if for any vertex $x$ of $C$, the set of faces 
incident with $x$ is a connected 
edge set of the dual graph of $\iota$.

\begin{thm}\label{embed_via_graph_vertical}
 Let $C$ be a simplicial complex. Then $C$ embeds into \Sthree\ if and only if its vertical split 
complex $\hat C$ has an embedding into \Sthree\ 
that 
satisfies the vertical dual graph connectivity constraints. 
\end{thm}

\begin{proof}
 First assume that $C$ has a topological embedding $\iota$ in \Sthree. Let $\iota'$ be the 
embedding  induced by $\iota$ of $\hat C$. 
By \autoref{same_complex}, $\iota$ and $\iota'$ have the same dual graph. Hence by 
\autoref{is_conny}, $\iota'$ satisfies the vertical dual graph 
connectivity constraints. 

Now conversely assume that $\iota'$ is an embedding into \Sthree\ of $\hat C$ that 
satisfies the vertical dual graph connectivity constraints. Let $G$ be the dual graph of $\iota'$. 
We shall 
recursively construct a sequence $(C_n)$ of 
simplicial complexes by identifying vertices that belong to the same vertex of $C$ that all have 
the 
vertical split 
complex $\hat C$ and topological embeddings $\iota_n$ of $C_n$ into \Sthree\ that all have the same 
dual graph $G$. 

If $C_n=C$, we stop and are done.
So there is a vertex $v$ of $C$ such that $C_n$ has at least two 
vertices 
cloned from $v$. The set of faces incident with $v$ is a connected edge set of $G$. So there are 
two distinct vertices $v_1$ and 
$v_2$ of $C_n$ cloned from $v$ whose incident faces 
share a vertex when considered as edge sets of 
$G$.
Hence there is a local surface of $\iota_n$ that contains $v_1$ and $v_2$. We obtain $C_{n+1}$ from 
$C_n$ by identifying $v_1$ and $v_2$.
The existence of a suitable embedding $\iota_{n+1}$ follows from \autoref{loc_surface_identify}. 

Since this recursion cannot continue forever, we must eventually have that $C_n=C$. Then $\iota_n$ 
is the desired embedding of $C$ and we 
are done.
\end{proof}

\subsection{Constructing embeddings from edge split complexes}

Our next step is to prove the following lemma analogously to one of the implications of  
\autoref{embed_via_graph_vertical}. 
Given a simplicial complex $C$ and a topological embedding $\iota$ into \Sthree\ of any of its 
split complex 
$\hat C$ into \Sthree, we say that 
$\iota$ satisfies the 
\emph{dual graph connectivity constraints} (with respect to $C$) if for any vertex or edge $x$ 
of $C$, the set of faces incident with 
$x$ is a connected edge 
set of the dual graph of $\iota$.

\begin{lem}\label{embed_via_graph_edge1}
Let $C$ be a locally connected simplicial complex. 
Assume that the split complex of $C$ has an embedding $\iota'$ into \Sthree\ 
that satisfies the dual graph connectivity constraints. 
 Then $C$ has an embedding in \Sthree\ that has the same dual graph as $\iota'$.
\end{lem}

Working with a strip instead of a unit interval, one shows the 
following analoguously to \autoref{loc_surface_identify}. 

\begin{lem}\label{loc_surface_identify2}
 Let $C$ be a simplicial complex obtained from a simplicial complex $C'$ by identifying two edges 
$e$ 
and $e'$ with disjoint sets of endvertices.
 Let $\iota'$ be a topological embedding of $C'$ into \Sthree. Assume that there is a local surface 
of $\iota'$ that contains both $e$ and 
$e'$. Then there is a topological embedding of $C$ into \Sthree\ that has the same dual graph as 
$\iota'$.
\qed
\end{lem}

\begin{proof}[Proof of \autoref{embed_via_graph_edge1}.]
Since the split complex is independent of the ordering in which we do splittings, the split 
complex $C'$ of $C$ is obtained by a sequence of the following operations:
first we split an edge. Then we split the two endvertices of that edge. After that the complex is 
again locally connected. So we eventually derive at the split complex.

We make an inductive argument similary as in the proof of \autoref{embed_via_graph_vertical}. Thus 
it suffices to show that if a complex embeds and satisfies the dual graph connectivity constraints 
at the clones of some edge, we can reverse the splitting at that edge within the embedding. 

After such a splitting operation the original edge is split into a set of vertex-disjoint edges. By 
the dual graph connectivity constraints, there are two of these edges in a common local surface of 
the embedding. So we can apply \autoref{loc_surface_identify2} to identify them. Arguing 
inductively, we can identify them all recursively. This shows why one such splitting can 
be reversed. Hence we can argue inductively as in the proof of \autoref{embed_via_graph_vertical} 
to complete the proof.
\end{proof}

\subsection{Embeddings induce embeddings of split complexes}

The goal of this subsection is to prove the following. 

\begin{lem}\label{embed_via_graph_edge}
Let $C$ be a locally connected simplicial complex with an embedding $\iota$ in \Sthree. 
Then its split complex has an embedding into \Sthree\ 
that satisfies the dual graph connectivity constraints and has the same dual graph as $\iota$. 
\end{lem}

Before we can prove this, we need some preparation. We start with the following lemma very similar 
to \autoref{loc_surface_identify2}. We define `detemined' and reveal the definition in the proof of 
the next lemma. 

\begin{lem}\label{loc_strip}
 Let $C$ be a simplicial complex obtained from a simplicial complex $C'$ by identifying two edges 
$e$ 
and $e'$ that only share the vertex $v$.
 Let $\iota'$ be a topological embedding of $C'$ into \Sthree. Assume that the embedding of $L(v)$ 
in the plane induced by $\iota'$ has a region\footnote{Component of $\Sbb^2$ without 
$L(v)$}  that contains both $e$ and 
$e'$. Then there is a topological embedding of $C$ into \Sthree\ that has the same dual graph as 
$\iota'$. The cyclic orientation at the new edge is determined. 
\end{lem}

\begin{proof}
We image that the link graph at $v$ is embedded in a small ball around $v$. Then the region $R$
containing $e$ 
and $e'$ is included in a unique local surface of $\iota'$. We call that local surface $\ell$.
We obtain $\bar C$ from $C'$ by adding a face $f$ at the edges $e$, $e'$ and one new edge.
The embedding $\iota$ induces an embedding of $\bar C$ as follows. We embed $C'$ as prescribed by 
$\iota'$ and embed $f$ in $\ell$. It remains to specify the faces just before or just after  $f$ 
at 
$e$ and $e'$. The face $f'$ just before $f$ at $e$ corresponds to some edge of $L(v)$ that has the 
region $R$ on its left, when directed towards $e$. Similarly, the face $f''$ just after $f$ at $e'$ 
corresponds to some edge of $L(v)$ that has the 
region $R$ on its right, when directed towards $e'$.
This embedding of $\bar C$ induces some embedding of $C$ by first contracting the third edge of 
$f$, the one not equal to $e$ or $e'$ and then contracting the face $f$, that is, we identify $e$ 
and $e'$ along $f$. Clearly this embedding has the same dual graph as $\iota'$. 

It remains to show that the cyclic orientation of the incident faces induced by the embedding 
at the new edge is determined. For that we reveal the 
definition of determined. It means that the cyclic ordering at the new edge is obtained by 
concatenating the cyclic orientations of $e$ and $e'$ induced by $\iota'$ so that $f'$ is followed 
by $f''$.  
\end{proof}

For the rest of this subsection we fix a topological embedding $\iota$ of a 
locally connected simplicial 
complex $C$ into \Sthree. 
Our aim is to explain how $\iota$ gives rise to an embedding of any split complex of $C$. First 
we need some preparation. 
 Let $\Sigma=(\sigma(e)|e\in E(C))$ be the combinatorial 
embedding induced by $\iota$.

Let $e$ be an edge of $C$ and $I$ a subinterval of $\sigma(e)$.
Let $\bar C$ be the simplicial complex obtained from $C$ by replacing $e$ by two edges, one that is 
incident with the faces in $I$ and the 
other that 
is incident with the faces incident with $e$ but not in $I$. 
We call $\bar C$ the simplicial complex obtained from $C$ by \emph{opening the edge $e$ along $I$}. 
We refer to the two new edges as the \emph{opening clones} of $e$. If we apply several openings, we 
extend the 
notion of opening cloning iteratively so that each edge of the resulting simplicial complex is 
opening cloned from a unique edge of $C$.

Let $C'$ be a simplicial complex obtained from a simplicial complex $C$ by splitting edges.
Given a rotation system $\Sigma$ of $C$, we obtain the \emph{induced} rotation system of $C'$ 
by restricting for each $e'$ of $C'$ 
cloned from an edge $e$ of $C$ the cyclic ordering $\sigma(e)$ to the faces incident with $e'$.
We define also an induced rotation system if $C'$ is obtained from $C$ by opening edges. This is 
as above with `clone' replaced by 
`opening clone'.

Let $\bar C$ be a simplicial complex obtained from $C$ by opening an edge and let $\bar\Sigma$ be 
the rotation system induced by 
$\Sigma$. 

\begin{lem}\label{open_edge}
The simplicial complex $\bar C$ has a topological embedding $\bar \iota$ into \Sthree\ whose 
induced planar rotation system is
$\bar \Sigma$.

The dual graph of $\bar \iota$ is obtained from the dual graph $G$ of $\iota$ by identifying 
the two endvertices of $I$ when considered as a trail in $G$.
\end{lem}

In particular, if $I$ is a closed trail in $G$, then $G$ is the dual graph of $\bar \iota$.
\begin{proof}[Proof of \autoref{open_edge}.]
 We can modify the embedding of $C$ such that there is an open cylinder around $e$ that does not 
intersect any edge except for $e$ or any face not incident with $e$. And all faces in $I$ intersect 
that cylinder only in the left half of the cylinder and the others only in the right half. Now we 
replace $e$ by 
two copies - one in the left half, the other in the right half. It is straightforward to check that 
the dual graph of the embedding has the desired property. 
\end{proof}

We fix an edge $e$ of $C$ with endvertices $v$ and $w$. 
\begin{lem}\label{make_interval}
There is an embedding $\iota'$ of $C$ in \Sthree\ that has the same dual graph as $\iota$ 
such that there is some connected component $X$ at $e$ that is a subinterval of the 
cyclic orientation $\sigma'(e)$, where $\Sigma'=(\sigma'(e)|e\in E(C))$ is the induced rotation 
system of $\iota'$. 
 \end{lem}
 
\begin{eg}
The following example demonstrates that in \autoref{make_interval} we cannot always pick 
$\iota'=\iota$. In the embedding in 3-space indicated in \autoref{fig:edge} no component 
at the edge $e$ is a subinterval of the cyclic orientation of the faces incident 
with $e$ induced by the embedding.

   \begin{figure} [htpb]   
\begin{center}
   	  \includegraphics[height=3cm]{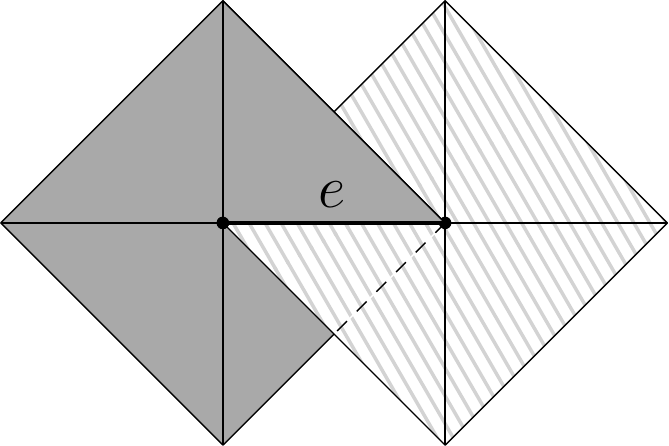}
   	  \caption{This complex is obtained by gluing together two discs, each with four faces, at 
the edge $e$.}\label{fig:edge}
\end{center}
   \end{figure}

\end{eg}
 Before we can prove \autoref{make_interval}, we need some preparation.

Given a cyclic orientation $\sigma$ and a subset $X$, we say that two elements $y_1$ and $y_2$ of 
$\sigma$ \emph{separate} $X$ for $\sigma$ if they are both not in $X$ and the two 
intervals\footnote{By $y_1\sigma y_2$ we denote the subinterval of $\sigma$ starting at $y_1$ and 
ending with $y_2$. } 
$y_1\sigma y_2$ and $y_2\sigma y_1$ both contain elements of $X$.

\begin{lem}\label{is_interval}
Let $\sigma$ be a cyclic orientation and $(P_i|i\in [n])$ be a partition of the 
elements of $\sigma$ such that no two elements of the same $P_i$ separate some other $P_j$.
Then there is some $P_k$ that is a subinterval of $\sigma$.
\end{lem}

\begin{proof}
We pick an arbitrary element $a$ of $P_1$. We may assume that a partition class $P_2$ exists.
For any $P_i$ not containing $a$, we define its 
\emph{first element} to first element of $P_i$ after $a$ in $\sigma$, and its \emph{last 
element} to first element of $P_i$ before $a$ in $\sigma$. 
The \emph{closure} of $P_i$ consists of those elements of $\sigma$ between its first and last 
element (including the first and the last one). We denote the closure of $P_i$ by $\overline{P_i}$. 

By assumption any two such closures $\overline{P_i}$ and $\overline{P_j}$ are either disjoint or 
contained in one 
another, that is, $\overline{P_i}\se \overline{P_j}$ or vice versa. 
Let $P_k$ be such that its closure is inclusion-wise minimal. Then $P_k$ is equal to its closure 
and hence a subinterval of $\sigma$.
\end{proof}

\begin{figure} [htpb]   
\begin{center}
   	  \includegraphics[height=3cm]{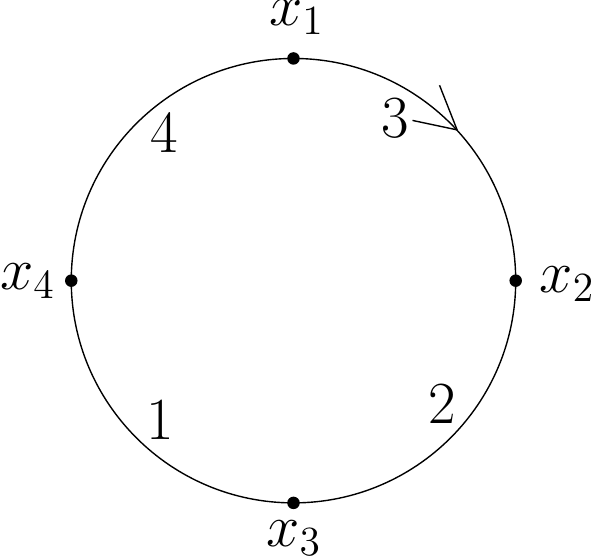}
   	  \caption{The cyclic orientation $\sigma$ is depicted as a cycle. The four segments 
between the element $x_i$ are labelled with the elements of $\Zbb_4$. This describes the ordering 
in which these segments are travered by the exchanged cyclic orientation.}\label{fig:sigma}
\end{center}
   \end{figure}

Given $e\in \sigma$, we denote the element just before $e$ 
by $e-1$ and the element just after $e$ by $e+1$. 
Given a cyclic orientation $\sigma$ and four of its elements $x_1$, $x_2$, $x_3$, $x_4$ such that 
$(x_1x_2x_3x_4)$ is a cyclic subordering of $\sigma$, the \emph{exchange} of $\sigma$ 
with respect to $x_1, x_2,x_3,x_4$ is the following cyclic orientation on the same 
elements as 
$\sigma$. We concatenate the two cyclic orientations obtained from $\sigma$ by deleting 
$x_1\sigma x_3-x_1-x_3$ and $x_3\sigma x_1$ such that the immediate successor of $x_4$ is $x_2$;
see \autoref{fig:sigma}, formally, it is \[
               x_3\sigma x_4(x_2\sigma x_3-x_3)(x_1\sigma x_2-x_1-x_2) (x_4+1)\sigma x_3  
                \]

Let $(P_i|i\in I)$ be a partition of the elements of $\sigma$, the \emph{fluctuation} of $\sigma$ 
with respect to $(P_i|i\in I)$ is the number of adjacent elements of $\sigma$ in different $P_i$.
Given a partition $\Pcal=(P_i|i\in I)$ of $\sigma$, an exchange is \emph{$\Pcal$-improving} 
if $x_2$ and $x_4$ are in the same $P_i$ but none of the following four pairs is in the same $P_i$:
$(x_4, x_4+1)$, $(x_2,x_2-1)$, $(x_1, x_1+1)$, $(x_3,x_3-1)$.

\begin{lem}\label{is_improving}
A cyclic orientation $\sigma'$ obtained from $\sigma$ by an exchange that is $\Pcal$-improving has 
strictly smaller fluctuation.
\end{lem}

\begin{proof}
The adjacent elements of $\sigma$ and $\sigma'$ are the same except for four pairs involving $x_1, 
x_2,x_3,x_4$. For $\sigma$ these pairs are those mentioned in the definition of 
`$\Pcal$-improving'. All these four pairs contribute to the fluctuation by the 
definition of $\Pcal$-improving. For $\sigma'$ the pair $(x_4,x_2)$ does not contribute to the 
fluctuation. 
\end{proof}

One way to partition $\sigma(e)$ is to put two elements of $\sigma(e)$ in the same class if -- when 
considered as edges of $L(v)$ -- they have endvertices in the same component of $L(v)-e$. An 
exchange is \emph{ $v$-improving} for $\sigma(e)$ if it is $\Pcal$-improving for that 
particular partition.

For the next lemma we fix the following notation. Let $X$ be the set of edges between $e$ and a 
connected component of $L(v)-e$.
Let $Y$ be the set of edges between $e$ and a connected component of $L(w)-e$.
Assume that no connected component at $e$ includes both $X$ and $Y$.

\begin{lem}\label{exchange_topo}Assume that two elements of $Y$ 
separate $X$ in the cyclic orientation $\sigma(e)$.
Then there is an embedding $\iota'$ of $C$ in \Sthree\ that has the same dual graph as $\iota$ 
such that $\sigma'(e)$ is obtained from $\sigma(e)$ by a $v$-improving exchange, where 
$\Sigma'=(\sigma'(e)|e\in E(C))$ is the induced rotation system 
of $\iota'$. 
\end{lem}

\begin{proof}
We abbreviate $\sigma(e)$ by $\sigma$. 
We denote the connected component at $e$ including $Y$ by $c(Y)$. 

\begin{sublem}\label{sublem99}
There are edges $f_1$ and $f_3$ of $c(Y)$ that separate $X$ such that the region of $L(w)$ just 
after $f_1$ is equal to the region just before $f_3$. And $f_1+1$ and $f_3-1$ are not in 
$c(Y)$.  
\end{sublem}

\begin{proof}
Let $f_1'$ and $f_3'$ be two elements of $Y$ that separate $X$. 
We fix two elements $x_1$ and $x_2$ of $X$ such that $x_1$ is in $f_1'\sigma f_3'$ and $x_2$ is in 
$f_3'\sigma  f_1'$. By choosing $f_1'$ and $f_3'$ as near to $x_1$ as possible, we ensure that the 
region just after $f_1'$ is equal to the region just before $f_3'$. We denote this region by $R$.
Let $Y'$ be the set of edges between $e$ and a connected component of $L(w)-e$ that is included in 
$c(Y)$. The set of all such $Y'$ is denoted by $\Ycal$. 
By replacing $Y$ by any $Y'\in \Ycal$ if necessary, we may assume that no set $Y'\in \Ycal$ 
contains 
elements both before and after $x_1$ on $f_1'\sigma f_3'$; indeed, by any such replacement 
$f_1'\sigma f_3'$ strictly decreases. 

\begin{sublem}\label{f1}
The interval  $f_1'\sigma x_1$ contains some $f_1\in c(Y)$ such that the 
region just after $f_1$ is $R$ and $f_1+1$ is not in $c(Y)$.
\end{sublem}

\begin{proof}
We recursively define a sequence $f_1^n$ of elements of  $f_1'\sigma x_1$. They are strictly 
increasing and contained in $c(Y)$. 
We start with $f_1^1=f_1'$. Assume that we already constructed $f_1^n$.
If $f_1^n+1$ is not in $c(Y)$ we stop and let $f_1=f_1^n$. Otherwise $f_1^n+1$ is in $c(Y)$. Let 
$Y'\in \Ycal$ so that $f_1^n+1\in Y'$. 

We prove inductively during this construction that any set $Y''\in \Ycal$ that contains an element 
of $(f_1^n+1)\sigma x_1$ contains no element of $f_1'\sigma f_1^n$. 

By the induction hypothesis, $Y'$ is a subset of $(f_1^n+1)\sigma x_1$. Let $f_1^{n+1}$ be the 
maximal element of $Y'$ in $(f_1^n+1)\sigma x_1$. By construction $f_1^{n+1}\in c(Y)$ and 
 $f_1^{n+1}$ is strictly larger than $f_1^n$. The region $R$ is just before $f_1^n+1$, the first 
element of $Y'$. Thus the region after $f_1^{n+1}$, the last element of $Y'$, must also be $R$. 
The induction step follows from the planarity of $L(w)$ as there is a component of $L(w)-e$ that is 
adjacent to the set $Y'$, and the induction hypothesis.

This process has to stop as $f_1'\sigma x_1$ is finite and the $f_1^n$ are strictly increasing. 
Thus we eventually find an $f_1$.
\end{proof}

Similarly as \autoref{f1} one shows that the interval $x_1\sigma f_3'$ contains some 
$f_3\in c(Y)$ such that the region just before $f_3$ is $R$ and $f_3-1$ is not in 
$c(Y)$.
So $f_1$ and $f_3$ have the desired properties. 
\end{proof}

We obtain $C_1$ from $C$ by opening the edge $e$ at the subinterval $f_1\sigma f_3$ of 
$\sigma$. By $\iota_1$ we denote the embedding of $C_1$ induced by $\iota$. 
By the choice of $f_1$ and $f_3$, the local surface just after $f_1$ is equal to the 
local surface just before $f_3$. Hence by  \autoref{open_edge} the embeddings $\iota_1$ and $\iota$ 
have the same dual graph. 

By \autoref{sublem99}, the link graph at $w$ of $C_1$ has two connected 
components. 
We obtain $C_2$ from $C_1$ by splitting the vertex $w$. 
By $\iota_2$ we denote the 
embedding of $C_2$ induced by $\iota_1$. As splitting vertices does 
not change the dual graph by \autoref{same_complex}, the embeddings $\iota_2$ and $\iota_1$ have 
the same dual graph. Summing up, $\iota_2$ and $\iota$ have the same dual graph.

We denote the copy of $e$ incident with $f_1$ by $e'$ and the other copy by $e''$. 
Since $e'$ and $e''$ are both incident with edges of $X$, the component of $L(v)-e$ adjacent to 
the  edges of $X$ has in the link graph of $C_2$ the two vertices $e'$ and $e''$ in the 
neighbourhood. Thus the vertices $e'$ and $e''$ share a 
face in the link graph at $v$ of $C_2$.

By \autoref{loc_strip} 
$\iota_2$ induces an embedding $\iota'$ of $C$ in \Sthree\ that has the same dual graph as 
$\iota_2$. 
Let $\Sigma'=(\sigma'(e)|e\in E(C))$ is the induced rotation 
system of $\iota'$. We denote the element of $X$ in $f_1\sigma f_3$ nearest to $f_1$ by $f_2$. 
Similarly, by $f_4$ we denote the element of $X$ in $f_3\sigma f_1$ nearest to $f_1$.
As $\sigma'(e)$ is determined by \autoref{loc_strip}, it is obtained by concatenating the cyclic 
orientations at $e'$ and $e''$ 
so that $f_4$ is followed by $f_2$. 
 That is, $\sigma'(e)$ is obtained from $\sigma(e)$ by exchanging with respect to $f_1, f_2, f_3, 
f_4$. 

It remains to check that this exchange is $v$-improving. Both $f_2$ and $f_4$ are in $X$. On the 
other hand $f_1$ and $f_3$ are in $c(Y)$ but $f_1+1$ and $f_3-1$ are not in $c(Y)$. In particular, 
they are in different $P_i$. Whilst $f_2$ and $f_4$ are in $X$, the two elements $f_2-1$ and 
$f_4+1$ are not in $X$. Thus this exchange is $v$-improving. 
\end{proof}

\begin{proof}[Proof of \autoref{is_interval}.]
By $(R_k|k\in K)$ we denote the 
partition of the faces incident with $e$ into the connected components at $e$. 
If no two elements of the same $R_a$ separate some other $R_b$, then by \autoref{is_interval} there 
is some $R_a$ that is a subinterval of $\sigma(e)$. In this case we can just pick $\iota'=\iota$ 
and are done.

We define the partition $(P_i|i\in I)$ of the faces incident with $e$ as follows. Two faces 
incident with $e$ are in the same partition if -- when considered as edges of $L(v)$ -- they have 
endvertices in the same component of $L(v)-e$. We define the partition $(Q_j|j\in J)$ the same with 
`$w$' in place of `$v$'. If some $P_i$ contains two elements separating some $Q_j$ for the 
cyclic orientation at $e$, we can apply \autoref{exchange_topo} to construct a new embedding 
of $C$. We do this until there are no longer such pairs $(P_i,Q_j)$. This has to stop after 
finitely many steps as by \autoref{is_improving} the fluctuation -- which is a non-negative 
constant only defined in terms of $(P_i|i\in I)$ -- of the cyclic orientation at $e$ strictly 
decreases in each step. So there is an embedding $\iota'$ of 
$C$ in \Sthree\ such that no  $P_i$ contains two elements separating some $Q_j$ for the 
cyclic orientation $\sigma'(e)$ and such that $\iota'$ has the same dual graph as $\iota$;
here we denote by $\Sigma'=(\sigma'(e)|e\in E(C))$ is the induced rotation 
system 
of $\iota'$.  
Hence by applying \autoref{is_interval}, it suffices to prove the following.
 
\begin{sublem}\label{R_topq}
For $\sigma'(e)$, either there is some $P_i$ containing two elements separating some $Q_j$ or no 
two elements of the same $R_a$ separate some other $R_b$.
\end{sublem}
 
 \begin{proof}
We assume that there is some $R_a$ that contains two elements $r_1$ and $r_2$ that separate some 
other $R_b$. The set $R_b$ is a disjoint union of sets $P_i$. Either $r_1$ and $r_2$ separate one 
of these $P_i$ or by the definition of connected component at $e$, there is some $Q_j$ included in 
$R_b$ that contains elements of different $P_i$, one included in $r_1\sigma'(e) r_2$ and the 
other in $r_2\sigma'(e) r_1$. Summing up there is some $P_i$ or $Q_j$ included in $R_b$ that is 
separated by $r_1$ and $r_2$. 

First we consider the case that there is a set $P_i$. So two 
elements of that set $P_i$ separate $R_a$. By an argument as above we conclude that there is some 
$P_m$ or $Q_n$ included in $R_a$ that is separated by two elements of $P_i$. 

Since the sets $P_m$ are defined from components of $L(v)-e$ and $\Sigma'$ induces an embedding of 
$L(v)$ in the plane, these components cannot attach at $e$ in a `crossing way', that no two 
elements of some $P_i$ can separate some other $P_m$. Thus there has to be such a set $Q_n$.

Summing up, if there is a set $P_i$ separated by $r_1$ and $r_2$, then it contains two elements 
separating some $Q_n$. Analogously one shows that otherwise the set $Q_j$ separated by $r_1$ and 
$r_2$ contains two elements 
separating some $P_n$. But then two elements of $P_n$ separate $Q_j$. This completes the proof. 
 \end{proof}

 By the construction of $\iota'$, no 
two elements of the same $R_a$ separate some other $R_b$ for $\sigma'(e)$. Then by 
\autoref{is_interval} there 
is some $R_a$ that is a subinterval of $\sigma'(e)$, as desired.  
\end{proof}

Let $C'$ be a simplicial complex obtained from the locally connected simplicial complex $C$ by 
splitting the edge $e$.  

\begin{lem}\label{split_once}
 There is a topological embedding $\iota'$ of $C'$ whose induced planar rotation system is the 
rotation system induced by $\Sigma$.
 
Moreover $\iota$ and $\iota'$ have the same dual graph.
\end{lem}

\begin{proof}

We denote the dual graph of $\iota$ by $G$. 
We prove this lemma by induction on the number of connected components at $e$. If there is only one 
such 
component, then $C'=C$ and the lemma is trivially true. So we may assume that there are at least 
two components.  By changing the embedding if necessary, by \autoref{make_interval} we may assume 
that there is a component $J$ at $e$ that is a subinterval 
of $\sigma(e)$. As $J$ is a subinterval of the closed trail 
$\sigma(e)$ of $G$, it is a trail in $G$. Next we show that it is a closed one:

 \begin{sublem}\label{is_closed}
The interval $J$ is a closed trail in $G$. 
\end{sublem}
\begin{proof}
We are to show that the local surface of the embedding just before the first face $f_1$ of $J$ is 
the 
same as the local surface just after the last edge $f_2$ of $J$. For that it suffices to show 
that in the embedding of the 
link graph $L(v)$ of $v$ induced by $\Sigma$, the region just before the edge $f_1$ is the same as 
the region just after the edge $f_2$. This 
follows from the fact 
that $J$ is the set of edges out of a set of connected components of 
$L(v)-e$. Indeed, the first and last edge out of every component are always in the same region. 
\end{proof}

We obtain $\bar C$ from $C$ by opening the edge $e$ along $J$. 
By \autoref{open_edge}, $\bar C$ has a topological embedding $\bar \iota$ into \Sthree\ whose 
induced planar rotation system is induced by $\Sigma$. By \autoref{is_closed}
and 
\autoref{open_edge}, the dual graph of $\bar \iota$ is $G$. 

We observe that $C'$ is obtained from $\bar 
C$ by splitting the clone of $e$ that corresponds to the subinterval $\sigma(e)\sm J$. Thus the 
lemma follows by applying induction on $\bar C$ and 
$\bar \iota$.
 \end{proof}

\begin{proof}[Proof of \autoref{embed_via_graph_edge}.]
 The split complex of $C$ is obtained from $C$ by a sequence of edge splittings and vertex 
splittings. By changing the order of the splittings if necessary, we may assume that the complex is 
always locally connected before we perform an edge splitting.  Hence we can apply 
\autoref{split_once} and \autoref{embed_via_graph_vertical} recursively to construct an embedding of 
the split complex. Since in each splitting step the dual graph is preserved, it satisfies the dual 
graph connectivity constraints by \autoref{is_conny} applied to the dual graph of $\iota$. 
\end{proof}

\subsection{Proof of \autoref{embed_via_matroid}}
 
We summarise the results of the earlier subsections in the following. 

\begin{thm}\label{embed_via_graph}
 Let $C$ be a simplicial complex and $\hat C$ be its split complex. Then $C$ embeds into 
\Sthree\ if and only if $\hat C$ 
has an embedding into \Sthree\ that 
satisfies the dual graph connectivity constraints. 
\end{thm}

\begin{proof}
Assume that $C$ embeds into 
\Sthree. Then by \autoref{embed_via_graph_vertical} its vertical split complex embeds into \Sthree\ 
and satisfies the vertical graph connectivity constraints. 
Since the vertical split complex is locally connected, we can apply \autoref{embed_via_graph_edge} 
to get the desired embedding of the split complex. Note that this embedding has the same dual graph 
as the vertical split complex. Hence it also satisfies the connectivity constraints for the 
vertices. 

Now conversely assume that the split complex has an embedding $\iota'$ that satisfies the dual 
graph connectivity constraints. By \autoref{embed_via_graph_edge1} the vertical split complex has 
an embedding in \Sthree. As this embedding has the same dual graph as $\iota'$, it satisfies the 
vertical dual graph connectivity constraints. So we can apply \autoref{embed_via_graph_vertical}. 
This completes the proof. 
\end{proof}

Now we show how \autoref{embed_via_graph} implies \autoref{embed_via_matroid}.

\begin{proof}[Proof of \autoref{embed_via_matroid}]
 Let $C$ be a globally 3-connected simplicial complex and let $\hat C$ be its split 
complex. 
  If $C$ embeds into \Sthree, then $\hat C$ has an embedding into \Sthree\ whose dual graph $G$ 
satisfies the dual graph connectivity constraints by \autoref{embed_via_graph}. 
By \autoref{same_matroid2}, the two simplicial complexes $C$ and $\hat C$ have the same dual 
matroid. So by \autoref{3dual_matroid} the cycle matroid of $G$ is the dual matroid of $C$.
This completes the proof of the `only if'-implication. 

Conversely assume that a split complex $\hat C$ of a simplicial complex $C$ has an embedding 
$\hat\iota$ 
into \Sthree\ 
and the dual matroid $M$ of $C$ is the cycle matroid of a graph $G$ and the set of faces incident 
with 
any vertex or edge of $C$ 
is a connected edge 
set of 
$G$. By \autoref{same_matroid2} $M$ is the dual matroid of $\hat C$. Let $G'$ be the dual graph of 
the embedding $\hat\iota$ of $\hat C$. By \autoref{3dual_matroid} the cycle matroid of $G'$ is 
equal to $M$. Since $M$ is 3-connected by assumption, by a theorem of Whitney \cite{Whitney_flip}, 
the graphs $G$ and $G'$ are identical. Hence $G'$ satisfies the connectivity constraints. So we can 
apply the 
`if'-implication of \autoref{embed_via_graph} to deduce the  `if'-implication of 
\autoref{embed_via_matroid}.
\end{proof}

\begin{proof}[Proof of \autoref{combined}.]
By \cite{3space1}, it suffices to show that a simplicial complex $C$ whose split 
complex is 
embeddable has an embedding if and only if its dual matroid has no constraint minor in the list 
of \autoref{fig:constraint}. Since the split 
complex is embeddable, its dual matroid is the cycle matroid of a graph $G$. By 
\autoref{same_matroid2} the dual matroid of $C$ is the cycle matroid of $G$.
By \autoref{embed_via_matroid}, $C$ is embeddable if and only if $G$ satisfies the graph 
connectivity constraints. The later is true if and only if there is no vertex or edge such that the 
set $X$ of incident faces is disconnected in $G$. 
By the main result of \cite{3space3}, $X$ is disconnected in $G$ if 
and only if $(G,X)$ has a 
constraint minor in the list of 
\autoref{fig:constraint}.
\end{proof}

\section{Infinitely many obstructions to embeddability into 3-space}\label{sec:examples}

In this section we construct an infinite sequence $(A_n|n\in \Nbb)$ of minimal obstructions to 
embeddability. More precisely, $A_n$ will have the property that its split complex is simply 
connected and embeddable, its dual matroid $M_n$ is the cycle matroid of a graph but no such 
graph will satisfy the connectivity constraints. However, if we remove a constraint or contract 
or delete an element from the dual matroid, then there is such a graph. 

\begin{figure} [htpb]   
\begin{center}
   	  \includegraphics[height=2cm]{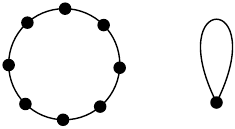}
   	  \caption{The matroid $M_8$. For each of the eight vertices on the cycle, there is a 
connectivity constraint forbidding that the loop is attached at that 
vertex.}\label{Mn}
\end{center}
   \end{figure}

The dual matroid $M_n$ of $A_n$ will be the disjoint union of a cycle $C_n$ of length $n$ and a 
loop $\ell$, see \autoref{Mn}. 

The connectivity constraints are as follows. Fix a cyclic 
orientation $\{e_i| i\in \Zbb_i\}$ of the edges on $C_n$. We have a connectivity constraint for 
every $i\in [n]$, namely that $X[i,n]=C_n-e_{i}-e_{i+1}+\ell$ is a connected set.   

\begin{fact}\label{not_met}
There is no graph whose cycle matroid is $M_n$ that meets all the connectivity constraints $X[i,n]$.
\end{fact}

\begin{proof}
By $\overline{C_n}$ we denote the graph that is a cycle of length $n$ whose edges have the cyclic 
ordering $\{e_i| i\in \Zbb_i\}$. It is straightforward to see that $\overline{C_n}$ is the unique 
graph 
whose cycle matroid is $C_n$ that meets all the connectivity constraints
$X[i,n]-e$. 

Now suppose for a contradiction that there is a graph $G$ whose cycle matroid is $M_n$ that meets 
all the connectivity constraints $X[i,n]$. Then $G$ is obtained from $\overline{C_n}$ by attaching 
a 
loop. Since each $X[i,n]$ contains $e$, we have to attach the loop at some vertex of 
$\overline{C_n}$. 
The connectivity constraint $X[i,n]$, however, forbids us to attach the loop at the vertex incident 
with $e_i$ and $e_{i+1}$. Hence $G$ does not exist. 
\end{proof}

A careful analysis of this proof yields the following simple facts.

\begin{fact}\label{met_rem}
\begin{enumerate}
\item 
There is a graph whose cycle matroid is $M_n$ that meets all the connectivity constraints $X[i,n]$ 
but one.
\item for every element $e$, there is a graph whose cycle matroid is $M_n-e$ that meets all the 
connectivity constraints $X[i,n]-e$;
\item for every element $e$, there is a graph whose cycle matroid is $M_n/e$ that meets all the 
connectivity constraints $X[i,n]-e$.
\end{enumerate}

\qed
\end{fact}

Hence it remains to construct $A_n$ such that its dual matroids is $M_n$ and so that the 
nontrivial connectivity constraints are the $X[i,n]$. We remark that we allow the faces of $A_n$ to 
be arbitrary closed walks. (One obtains a simplicial complexes from $A_n$ by applying baricentric 
subdivisions to the faces.) 

\begin{figure} [htpb]   
\begin{center}
   	  \includegraphics[height=3cm]{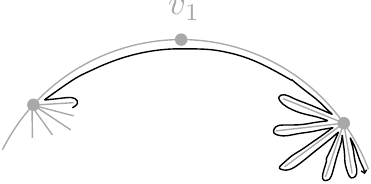}
   	  \caption{In grey, we indicate the cycle $C$ with the new 
edges. In black we sketched the traversal of the face $e_1$ after addition of the new 
edges. }\label{e_1new}
\end{center}
   \end{figure}

We start the construction of $A_n$ with a cycle $C$ of length $n$. 
We attach $n$ faces, which we call $e_1,...,e_n$. For each $e_i$, and each vertex $v_k$ of $C$ 
except for the $i$-th vertex $v_i$, we attach $n-1$ edges and let $e_i$ traverse them in between 
the 
two edges incident with $v_k$. We denote the endvertices of the new edges not on $C$ by $x(i,k,j)$ 
where $(k,j\leq n; k,j\neq i)$, see \autoref{e_1new}. 

Next we disjointly add a copy of the original cycle $C$ and only attach a single face to it which 
we 
denote by $\ell$. Call the resulting walk-complex\footnote{A \emph{walk-complex} is a graph 
together with a family of closed walks, which we call its faces. Every simplicial complex is a 
walk-complex. Conversely, from every walk complex we can build a simplicial complex by attaching 
at each face a cone over that walk. } $A_n'$. We finally obtain $A_n$ from $A_n'$ 
by 
identifying for each $i\in [n]$ the $i-th$ vertex $v_i$ on the new copy of $C$ with all vertices 
$x(i',i,i)$ with $i'\neq i$. 

By construction, the split complex of $A_n$ is $A_n'$. Hence by \autoref{same_matroid2} above, 
the dual matroid of $A_n$ is $M_n$. By construction, the nontrivial connectivity constraints are 
the $X[i,n]$. Clearly, the split complex $A_n'$ is simply connected and embeddable. 

This completes the construction of the $A_n$. By \autoref{not_met} and \autoref{met_rem} they have 
the desired properties.

\appendix

\section{Appendix {I}}\label{appendixA}

First we give a definition of `globally 3-connected' directly in terms of the simplicial complex 
without 
referring to its dual matroid. 
Given a simplicial complex $C$, its edge/face incidence matrix $A$ and a subset $L$ of the faces of 
$C$, we 
denote by $r(L)$ the rank over $\Fbb_3$ of the submatrix of $A$ induced by the vectors whose faces 
are in $L$. 
A \emph{2-separation} of a simplicial complex $C$ is a partition of its set $F$ of faces into two 
sets $L$ 
and $R$ both of size at least two such that $r(L)+r(R)\leq r(F)+1$. It is straightforward it check 
that a simplicial complex is globally 3-connected if and only if it has no 2-separation.

\vspace{.3 cm}

When defining `edge split complexes', we mentioned a related more naive definition. Here we give 
this definition. In \autoref{several_edge_split_complexes} 
and \autoref{not_thm} we show that this notion lacks two important features of edge split complexes.
\emph{Splitting an edge $e$ at an endvertex $v$} is defined like `splitting $e$' but with 
`in the same connected component at $e$' replaced by `$v$-related'. A \emph{lazy edge split 
complex} is defined as `edge split complex' but with `for every edge there is only one 
component at $e$' replaced by `it is locally 2-connected'. \emph{lazy split complex} is defined 
like `split complex' with `lazy edge split complex' in place of `edge split complex'.

\begin{eg}\label{several_edge_split_complexes}
In this example we construct a simplicial complex $C$ that has two distinct lazy edge split 
complexes.
We will construct $C$ such that it has two vertices $v$ and $w$; these vertices are joined by 
five edges $e$, $e_1$, $e_2$, $e_3$ and $e_4$. The edge $e$ is a cut 
vertex in the link graphs at $v$ and $w$. And splitting $e$ at one endvertex will make the link 
graph at the other endvertex 2-connected, see \autoref{fig:not_unique}. 

   \begin{figure} [htpb]   
\begin{center}
   	  \includegraphics[height=3cm]{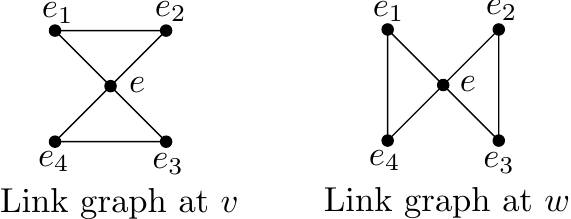}
   	  \caption{If we split one of these link graphs at $e$, the other becomes a six-cycle. 
}\label{fig:not_unique}
\end{center}
   \end{figure}

Next we construct $C$ with the above properties.  
We obtained $C$ from four triangular faces $f_1$, $f_2$, $f_3$ and $f_4$ glued together at a single 
edge $e$. Let $v$ and $w$ be the two 
endvertices of that edge.
Let $e_i[v]$ be the edge of $f_i$ incident with $v$ different from $e$. 
Let $e_i[w]$ be the edge of $f_i$ incident with $w$ different from $e$. 
Let $v_i$ be the vertex incident with $f_i$ that is not incident with $e$.
We add the edges $e_k$ between $v_k$ and $v_{k+1}$ for any $k\in \Zbb_4$. 
We add the four faces: $e_1[v]e_1e_2[v]$, $e_3[v]e_3e_4[v]$, $e_2[w]e_3e_3[w]$ and 
$e_4[w]e_4e_1[w]$. This completes the construction of $C$. 
\end{eg}

\begin{eg}\label{not_thm}
 In this example we show that \autoref{embed_via_graph} with `split complex' replaced by `lasy 
split complex' is false. Let $H$ be a planar graph with vertices $v$ and $w$ such that the graph 
$H'$ obtained from $H$ by identifying the vertices $v$ and $w$ is not planar. Let $C$ be the cone 
over $H$. We obtain $C'$ from $C$ by identifying the two edges corresponding to $v$ and $w$. Whilst 
the link at the top of $C$ is $H$, the complex $C'$ has the link $H'$ and is hence not embeddable. 
By choosing $v$ and $w$ far apart in $H$, one ensures that $C'$ is a simplicial complex. 

The lazy split complex of $C'$ is unique and equal to $C$. Unlike $C'$, the simplicial complex 
$C$ is embeddable. The dual graph of every embedding consists of a single vertex, and so 
trivially satisfies the graph connectivity constraints. This completes the example. 

Concerning \autoref{embed_via_matroid}, it is 
straightforward to modify the example to make the dual graph of the embedding 3-connected. 
\end{eg}

\section{Appendix {II}: Matrices representing matroids over the integers}\label{appendixB}

Matroids representable over the integers are well-studied \cite{oxley2}. In this appendix, 
we study something very related but slightly different, namely matrices that represent matroids 
over the integers. Our aim in this appendix is to prove \autoref{char_reg_represent} below, 
which is a characterisation of certain matrices representing matroids over the integers. 

A matrix $A$ is a \emph{representation of a matroid $M$ over a field $k$} if all its entries are in 
$k$ and the columns are indexed with the 
elements of $M$.
Furthermore for every circuit $o$ of $M$ there is a vector $v_o$ in the span over $k$ of the rows of 
$A$ whose support is $o$. 
And the vectors $v_o$ span over $k$ all row vectors of $A$.

\begin{comment}
 %Under these circumstances it follows that for every cocircuit $b$ of $M$, there is a vector $w_b$ 
that is orthogonal over $k$ to all row 
vectors of $A$ whose %support is equal to $b$. These vectors generate over $k$ all vectors that are 
orthogonal over $k$ to every row vector. 
\end{comment}

The following is well-known. 

\begin{lem}\label{contr_matrix}
 Let $A$ be a matrix representing a matroid $M$ over some field $k$. Let $I$ an element set that is 
independent in $M$.
Then the matrix obtained from $A$ by deleting all columns belonging to elements of $I$ represents 
the matroid $M/I$ over $k$. \qed
\end{lem}

A matrix $A$ is a \emph{regular representation} (or \emph{representation over the integers}) of a 
matroid $M$ if all its entries are integers and the columns are indexed with the 
elements of $M$.
Furthermore for every circuit $o$ of $M$ there is a $\{0,-1,+1\}$-valued vector\footnote{A 
\emph{vector} is an element of a vector space 
$k^S$, where $k$ is a field and $S$ is a set. In a slight abuse of notation, in this paper we also 
call elements of modules of the form 
$\Zbb^S$ vectors.} $v_o$ in the span over $\Zbb$ of the rows of $A$ whose support is $o$. 
And the vectors $v_o$ span over $\Zbb$ all row vectors of $A$.
The following is well-known. 

\begin{lem}\label{reg_repre_ortho}
Assume that a matrix $A$ regularly represents a matroid $M$.
Then for every cocircuit $d$ of $M$, there is a $\{0,-1,+1\}$-valued vector $w_d$ whose support is 
equal to $d$ that is 
orthogonal\footnote{Two vectors $a$ and $b$ in $k^S$ are \emph{orthogonal} if $\sum_{s\in S} 
a(s)\cdot b(s)$ is identically zero over $k$.} 
over $\Zbb$ to all row vectors of $A$. These vectors $w_d$ generate over $\Zbb$ all vectors that are 
orthogonal over $\Zbb$ to every row 
vector.\qed
\end{lem}
\begin{comment}
 We only add a proof sketch of the above well-known lemma. 

\begin{proof}[Proof sketch.] 
We fix a base $b$ of $M$.
First we show the fact that the family of vectors $v_o$, where $o$ is a fundamental circuit of $b$, 
generate over $\Zbb$ the  
whole span of the row vectors of $A$. To see this, note that we can subtract from any vector $a$ in 
the span a suitable $\Zbb$-valued 
combination of vectors $v_o$ of our family so that the remaining vector $a'$ is zero outside $b$. 
Since the support of $a'$ if nonempty 
includes a circuit of $M$, the support of $a'$ is empty and hence $a$ is in the span of the vectors 
in the family.  

For every fundamental cocircuit $d$ of $b$, let $z$ be the unique element of $d\cap b$.
We define the vector $w_d$ as follows. $w_d(z)=1$ and $w_d(x)$ is zero if $x\notin d$.
For $x\in d-z$, we define $w_d(x)$ so that $w_d$ gets orthogonal to the vector $v_o$ of the 
fundamental circuit $o$ of $b$ with respect to 
the element $x$, that is, $w_d(x)= - v_o(x) \cdot w_d(z) \cdot v_o(z)$. 
By construction the vectors $w_d$ are orthogonal to all vectors of the fundamental circuits of $b$. 
Hence by the above, they are orthogonal 
to all vectors in the span of $A$. 
An argument dual to above proof of the face shows that every vector orthogonal to the span of the 
row vectors of $A$ is spanned over $\Zbb$ 
by the $w_d$. 
\end{proof}
\end{comment}

The following is well-known. 

\begin{lem}\label{supp_contains_circuit}
 Let $M$ be a matroid regularly represented by a matrix $A$.
Let $v$ be a sum of row vectors of $A$ with integer coefficients.
If the support of $v$ is nonempty, then it includes a circuit of $M$.
\qed
\end{lem}

\begin{eg}
 A matrix is \emph{unimodular} if it is $\{0,-1,+1\}$-valued and the determinant of every quadratic 
submatrix is 
$\{0,-1,+1\}$-valued\footnote{Here we evaluate the determinate over $\Zbb$}.
Every unimodular matrix is a regular representation of some matroid, see for example 
\cite{Truemper_MatroidDecompositions}. 
For example, the vertex$\slash$edge incidence matrix of a graph $G$ is a regular representation of 
the graphic matroid of $G$. \begin{comment}
          Matroid Decompositions Lemma 9.2.1                   
                            \end{comment}
\end{eg}

 There also exist regular representations that are not totally unimodular:

\begin{eg}

\[\left(\begin{matrix}
 1  & 1\cr
1 & -1 \cr
1 & 0 \cr
0 & -1 \cr
\end{matrix}\right)
\]
This matrix is a regular representation of the matroid consisting of two elements in parallel but 
it is 
not totally unimodular. 
\end{eg}

A matroid is \emph{regular} if it can be regularly represented by some matrix. The class of regular 
matroids has many equivalent characterisations \cite{oxley2}.
For example, a matroid has a regular representation (in fact a 
totally unimodular one) if and only if 
it has a representation over every field.
In this paper, we need the following related fact, which focuses on the matrices instead of the 
matroids:

\begin{thm}\label{char_reg_represent}
  Let $A$ be a matrix whose entries are $-1$, $+1$ or $0$.
Then $A$ regularly represents a matroid if and only if there is a single matroid $M$ such that $A$ 
represents $M$ over any field. 
\end{thm}

Whilst the 'only if'-implication is immediate, the other implication is less obvious. 
To prove it we rely on the following.

\begin{lem}\label{loc_global}
 Let $(v_i|i\in I)$ be a family of integer valued vectors of $\Zbb^S$, where $S$ is a finite set.
Assume that the family $(v_i|i\in I)$ considered as vectors of the vector space $\Qbb^S$ spans the 
whole of $\Qbb^S$ over $\Qbb$.
Additionally, assume that for every prime number $p$, the same assumption is true with the finite 
field `$\Fbb_p$' in place of `$\Qbb$'.
Then the family $(v_i|i\in I)$ spans over $\Zbb$ all integer valued vectors in $\Zbb^S$.
\end{lem}

\begin{proof}[Proof that \autoref{loc_global} implies \autoref{char_reg_represent}.] 
Assume that $A$ is an integer valued matrix that represents the matroid $M$ over $\Qbb$ and over all 
finite fields $\Fbb_p$ for every prime 
number $p$, when we interpret\footnote{Here in $\Fbb_p$ we interpret the integer $m$ as its 
remainder after division by $p$.} the entries of 
$A$ as elements of the appropriate field. 
Our aim is to show that $A$ regularly represents the matroid $M$.

Let $b$ be a base of $M$. Let $A'$ be the matrix obtained from $A$ by deleting all columns belonging 
to elements of $b$. 
We denote by $M'$ the matroid $M/b$, in which every element is a loop. 
By \autoref{contr_matrix}, $A'$ represents the matroid $M'$ over $\Qbb$ and over all finite fields 
$\Fbb_p$. 
Let $(v_i|i\in I)$ be the family of row vectors of $A'$. 
Since every element of $M'$ is a loop, we can apply \autoref{loc_global} and deduce that the family 
$(v_i|i\in I)$ spans over $\Zbb$ all 
integer valued vectors in $\Zbb^{E'}$, where $E'$ is the set of elements of $M'$. 

Let $v$ be any integer valued vector that is generated by the rows of $A$ over $\Qbb$. 
We show that $v$ is also generated  by the rows of $A$ with integer coefficients. 
By the above, there is a vector $w$ generated from the row vectors of $A$ over $\Zbb$ that agrees 
with $v$ in all coordinates of $E'$. Hence 
$v-w$ is generated by the row vectors over $\Qbb$. So if $v-w$ is nonzero, its support must contain 
a circuit of $M$ by 
\autoref{supp_contains_circuit}. Since the support of $v-w$ is contained in the base $b$, the 
support does not contains a circuit of $M$. 
Hence $v$ must be equal to $w$. Thus $v$ is in the span of the row vectors with coefficients in 
$\Zbb$.
 
Now let $o$ be a circuit of $M$. Since $A$ is a regular representation of $M$ over $\Qbb$, there is 
a vector $v_o$ with entries in $\Qbb$ 
generated by the row vectors of $A$ over $\Qbb$ whose support is $o$. We multiplying all entries 
with a suitable rational number if 
necessary, we may assume that additionally
all entries of $v_o$ are integers and that the greatest common divisor of the entries is one. 
By the above $v_o$ is in the span of the row vectors with coefficients in $\Zbb$.

Next we show that all entries of $v_o$ are zero, plus one, or minus one. 
Suppose for a contradiction that there is some prime number $p$ that divides some entry of $v_o$. 
If we interpret the entries of $v_o$ as elements of $\Fbb_p$, then $v_o$ is also in the span of the 
row vectors with coefficients in 
$\Fbb_p$.
Indeed, the coefficients are just the integer coefficients we have in the representation over $\Zbb$ 
interpreted as elements of $\Fbb_p$.
Since the greatest common divisor of the entries of $v_o$ is one, $v_o$ when interpreted over 
$\Fbb_p$ is nonzero but its support is 
properly contained in $o$.
Since in $M$ the circuit $o$ does not include another circuit, we get a contraction to the 
assumption that $A$ represents $M$ over $\Fbb_p$.
Thus all entries of $v_o$ are zero, plus one, or minus one. 

It remains to show that the set of vectors $v_o$ where $o$ is a fundamental circuit of $b$ generates 
every row vector $x$ of $A$. 
Since for every element not in $b$, there is a unique $v_o$ which takes the value plus one or minus 
one 
at that element and zero at every other 
elements not in $b$, 
there is a vector $x'$ generated over $\Zbb$ by the $v_o$ that agrees with $x$ when restricted to 
$E'$. As above we deduce that $x'=x$, and 
hence $x$ is generated by the $v_o$ over $\Zbb$. 
Thus $A$ regularly represents $M$. 
\end{proof}

In order to prove \autoref{loc_global}, we rely on the following well-known lemma.
\begin{lem}\label{number_theory_lem}
 Let $m$ and $n$ be integer and let $d$ be their greatest common divisor. 
Then there are integers $\alpha$ and $\beta$ such that $\alpha\cdot m - \beta\cdot n=d$.
\qed
\end{lem}

\begin{proof}[Proof of \autoref{loc_global}.]
 Let $s\in S$ be arbitrary. By $e_s$ we denote the vector which in coordinate $s$ has the entry one 
and otherwise the entry zero. 
Since the family $(v_i|i\in I)$ spans $e_s$ over $\Qbb$, there is some positive natural number 
$\gamma_s$ so that the family $(v_i|i\in I)$ 
spans $\gamma_s\cdot e_s$ over $\Zbb$. 
Let $\delta_s$ be the least possible value for $\gamma_s$.
Our aim is to show that all $\delta_s$ are equal to one. 
Suppose not for a contradiction. Then there is some prime number $p$ that divides some $\delta_s$.
Let $\bar s$ be the index so that in the factorisation of $\delta_{\bar s}$ the prime number $p$ has 
the highest multiplicity, say $k$.

\begin{sublem}\label{sublem71}
There is some nonzero integer $\epsilon$ such that $p$ has the multiplicity at most $k-1$ in the 
factorisation of $\epsilon$ and such that 
$\epsilon\cdot e_{\bar s}$ is 
spanned by  the family $(v_i|i\in I)$ over $\Zbb$. 
\end{sublem}

Let us first see how we finish the proof assuming \autoref{sublem71}. By 
\autoref{number_theory_lem}, there are $\alpha$ and $\beta$ such 
that
$\alpha\cdot \delta_{\bar s} - \beta\cdot \epsilon$ is equal to the greatest common divisor $D$ of 
$\delta_{\bar s}$ and $\epsilon$.
Hence by \autoref{sublem71} $D\cdot e_{\bar s}$ is generated by the family $(v_i|i\in I)$ over 
$\Zbb$. Since $p$ has the multiplicity at 
most $k-1$ in the factorisation of $D$, the number $D$ is strictly smaller than $\delta_{\bar s}$. 
This contradicts the choice of 
$\delta_{\bar s}$.
Hence all $\delta_s$ are equal to one. 
It remains so show that the following.
\begin{proof}[Proof of \autoref{sublem71}.]
Since the family $(v_i|i\in I)$ spans $e_{\bar s}$ over $\Fbb_p$, there is an integer valued vector 
$w$ such that
the family $(v_i|i\in I)$ spans $e_{\bar s}+p\cdot w$ over $\Zbb$. 
For a subset $T$ of $S$ we denote by $w_T$ the vector which takes the value $w(s)$ in coordinate $s$ 
if $s\in T$ and zero otherwise. 
We denote the multiplicity of $p$ in the factorisation of an integer $n$ by $\sharp_p(n)$.

We shall show inductively for every subset $T$ of $S$ that there is some nonzero natural number 
$\epsilon_T$ with $\sharp_p(\epsilon_T)\leq 
k-1$ such that $\epsilon_T\cdot \left(e_{\bar s} +p\cdot w_T\right)$ is spanned by  the family 
$(v_i|i\in I)$ over $\Zbb$. 
We start the induction with $T=S$ and $\epsilon_T=1$ and so $w_T=w$. 
Assume that we already proved the induction hypothesis for a nonempty subset $T$ of $S$. Let $t\in 
T$ be arbitrary. 
We denote the greatest common divisor of $\epsilon_T\cdot p \cdot w(t)$ and $\delta_t$ by $d_t$.
We let $\epsilon_{T-t}= \epsilon_T\cdot \frac{\delta_t}{d_t}$. 
We have
\[
 \sharp_p(\epsilon_{T-t})= \sharp_p(\epsilon_T)+ \sharp_p(\delta_t)-\sharp_p(d_t)\leq 
\sharp_p(\epsilon_T)+ 
\sharp_p(\delta_t)-min\{\sharp_p(\epsilon_T)+1, \sharp_p(\delta_t)\}=
\]
\[
 =max\{\sharp_p(\delta_t),\sharp_p(\epsilon_T)-1\}
\]

Hence by the choice of $\bar t$ and by induction $\sharp_p(\epsilon_{T-t})\leq k-1$. 
Furthermore:
\[
  \frac{\delta_t}{d_t} \cdot \epsilon_T\cdot \left(e_{\bar t} +p\cdot w_T\right) - 
\frac{\epsilon_T\cdot p \cdot w(t)}{d_t}\cdot \delta_t 
e_t= \epsilon_{T-t}\cdot \left(e_{\bar t} +p\cdot w_{T-t}\right)
\]
Note that all fractions in the above equation are integers.
This completes the induction step.
Hence the vector $\epsilon_{\emptyset}\cdot e_{\bar t}$ is 
spanned by  the family $(v_i|i\in I)$ over $\Zbb$, which completes the proof. 
\end{proof}
\end{proof}

\bibliographystyle{plain}
\bibliography{literatur}

\end{document}